\input amstex.tex
\input amsppt.sty   
\magnification 1200
\vsize = 8.5 true in
\hsize=6.2 true in
\NoRunningHeads        
\parskip=\medskipamount
        \lineskip=2pt\baselineskip=18pt\lineskiplimit=0pt
       
        \TagsOnRight
        \NoBlackBoxes
        \def\ve{\varepsilon}
        
        \def\vp{\varphi}

        \topmatter
        \title
        Quasi Periodic Solutions of Nonlinear Random Schr\"odinger Equations
        \endtitle
\author
        J.~Bourgain\footnote"{$\dagger$}" {Partially supported by NSF grant DMS-03-22370\hfill\hfill\qquad}
 and W-M.~Wang
        \footnote"{$^{*}$}"{Partially supported by NSF grant
DMS-97-29992 and DMS 0503563\hfill\hfill\qquad}
        \endauthor        
\address{Institute for Advanced Study, Einstein Drive,
Princeton, N.J. 08540, U.S.A.}
\endaddress
\email
{bourgain\@math.ias.edu}
\endemail
\address
{D\'epartement de Math\'ematique, Universit\'e Paris Sud, 91405 Orsay Cedex, FRANCE}
\endaddress
        \email
{wei-min.wang\@math.u-psud.fr}
\endemail
\thanks{W.-M. Wang thanks the Institute for Advanced Study for hospitality, where part 
of this work was done.}
\endthanks

        \bigskip\bigskip
        \bigskip
        \toc
        \bigskip
        \bigskip 
        \widestnumber\head {Table of Contents}
      
        \head 1. Introduction and statement of the theorem\endhead
        \head 2. Hamiltonian representation and Lyapunov-Schmidt decomposition 
        \endhead
        \head 3. Newton scheme
        \endhead
        \head 4. $P$-equations and statements in $\theta$
        \endhead
        \head 5. Invertibility of $T(y_i)$, $Q$-equations and determination of $\omega$
        \endhead
        \head 6. Construction of $y_i$ and completion of the assemblage
        \endhead
        \head 7. Proof of the theorem
        \endhead 
	       \head 8. Appendix: Localization results for random Schr\"odinger operators
        \endhead 
        \endtoc
        \endtopmatter
        \vfill\eject
        \bigskip
\document
\head{\bf 1. Introduction and statement of the theorem}\endhead
\noindent {\it The nonlinear random Schr\"odinger equation.}

We seek time quasi-periodic solutions to the nonlinear random Schr\"odinger equation
$$
i\frac\partial{\partial t}u =(\epsilon\Delta+V)u+\delta |u|^{2p}u\qquad (p>0),\tag 1.1
$$
on $\Bbb Z^d\times [0,\infty)$, 
where $0<\epsilon, \delta\ll 1$, $\Delta$ is the discrete Laplacian:
$$\aligned \Delta_{ij}&=1,\quad |i-j|_{\ell^1}=1,\\
&=0,\quad \text{otherwise},\endaligned\tag 1.2$$
$V=\{v_j\}_{j\in\Bbb Z^d}$, the potential, is a family of {\it time independent} independently
identically distributed (i.i.d.)
random variables with common
distribution $g=\tilde g(v_j)dv_j$, $\tilde g\in L^{\infty}$.  We also assume $\text{supp }g$
is a bounded set. The probability space is taken
to be 
$$\Bbb R^{\Bbb Z^d} \text {with measure} \prod_{j\in\Bbb Z^d}g(v_j)=\prod_{j\in\Bbb Z^d}\tilde
g(v_j)dv_j,\,\tilde g\in L^{\infty}, \text{ supp }g \text{ is a bounded set}.\tag 1.3$$ 
$V=\{v_j\}_{j\in\Bbb Z^d}$ serve as parameters for the nonlinear problem in (1.1).

Given an initial condition $u(0)$ in $\ell^2(\Bbb Z^d)$, one of the central questions is whether
$u(t)$ remains localized for all $t$, i.e., if $u(0)\in \ell^2(\Bbb Z^d)$, $\forall \kappa$,
can one find $R$, such that
$$\Vert u(t)\Vert_{\ell^2(\{\Bbb Z\backslash[-R,R]\}^d)}<\kappa, \,\forall t?\tag 1.4$$
(From now on, we write $|\,|$ for $|\,|_{\ell^1}$, $\Vert\,\Vert$ for $\Vert\,\Vert_{\ell^2}$.) 
When $\epsilon=\delta=0$, the answer to (1.4) is affirmative. Since $u(0)=\sum_{j\in\Bbb
Z^d}a_j\delta_j$, $a_j\to 0$, as $|j|\to\infty$, $u(t)=\sum_{j\in\Bbb
Z^d}a_j\delta_je^{-iv_jt}$ is almost-periodic (infinite number of frequencies) and the 
upper bound in (1.4) is trivially verified.

In this paper, for appropriate initial conditions $u(0)$, we construct time quasi-periodic
solutions to (1.1). So the answer to (1.4) is affirmative for such $u(0)$'s. This is the
content of the Theorem and its Corollary. 

Before we enter into the heart of the matter, we first address question (1.4) to
\smallskip
\noindent {\it The linear random Schr\"odinger equation.}

When $\delta=0$, (1.1) reduces to the linear random Schr\"odinger equation:
$$\aligned
i\frac\partial{\partial t}u &=(\epsilon\Delta+V)u,\\
&{\overset\text{def }\to=}Hu\endaligned\tag 1.5
$$
on $\Bbb Z^d\times [0,\infty)$. When $0<\epsilon\ll 1$, it is well known from the works in 
\cite{AFHS, AM, vDK, FMSS, FS, GB, GK, GMP} etc. that the upper bound in (1.4) is
satisfied. This is customarily called Anderson localization (A.L.) after the physicist P. Anderson
\cite{An}. Since the potential is time independent: $V(j,t)=V(j)$, properties of time evolution can be
deduced from the spectral properties of $H$, which we summarize below. For more details, see the 
Appendix.

Let $\sigma(H)$ be the spectrum of $H$. For $H$ defined in (1.5), 
$$\sigma(H)=[-2\epsilon d, 2\epsilon d]+\text {supp}\, g,\, a.s.\tag 1.6$$
(Recall the probability space defined in (1.3).) \cite{CFKS, PF}.
If $0<\epsilon\ll 1$ and the probability measure satisfies (1.3), then almost surely the 
spectrum of $H$ is (dense) pure point, $\sigma(H)=\sigma_{pp}$, with exponentially localized
eigenfunctions: $\phi_j$, $j\in\Bbb Z^d$. 

Given $u(0)\in\ell^2(\Bbb Z^d)$, we decompose $u(0)$ as $u(0)=\sum_{j\in\Bbb
Z^d}a_j\phi_j$. So 
$$u(t)=\sum_{j\in\Bbb Z^d}a_j\phi_je^{-i\omega_jt},\tag 1.7$$ 
where $\omega_j$ are the eigenvalues for the eigenfunctions $\phi_j$. $u(t)$ is almost-periodic
and verifies the upper bound in (1.4). So equation (1.5) has A.L.
\smallskip
\noindent {\it Some motivations for studying equation (1.1).}

Schr\"odinger equations are equations that describe physical systems, which typically correspond
to a $n$-body problem. The linear equation in (1.5) is a $0^{\text {th}}$ order appproximation,
where the $n$-body interaction is lumped into the effective potential $V$. Quantum mechanically,
$|u|^2$ is interpreted as particle density, so the nonlinear term in (1.1) can be interpreted
as modelling particle-particle interaction. (The nonlinear term in (1.1) can be more general
and of convolution type. It will not affect our construction below.) This is sometimes called
the Hartree-Fock approximation (cf. \cite{LL, O, Sh}) and is a first order approximation to
the original $n$-body problem. This is our first motivation to study (1.1).
Other physical motivations along this line appear in \cite{FSW}.

In particular, our method permits us to construct quasi-periodic solutions for
the Landau-Lifschitz equations on
nonlinear classical spin waves with a large random external magnetic field.
Thus
$$
\dot S_j=S_j\times [(\Delta S)_j+h_j] \qquad (j\in\Bbb Z^d)
$$
where $S_j$ are unit vectors in $\Bbb R^3$ and $h_j=V_j\overset\rightarrow \to
e_3$ say;
with $V=(V_j)_{j\in\Bbb Z^d}$ a large random potential.

As explained in \cite{FSW}, we may then seek for a solution $S_j\approx e_3$ and
the perturbation is subject
to an equation of the form (1.1), but with a nearest neighbor convolution
nonlinearity instead of the local one
$|u|^{2p} u$ (see \cite{FSW} for details).
As mentioned before, (1.1) was chosen as a model but the method described in
the paper is sufficiently robust to
cover in particular any nonlinearity with finite range interactions.

Our second motivation originates from KAM type of stability questions for infinite dimensional
dynamical systems. (For results in the standard KAM context, see e. g. \cite{E}.) (1.1) is a Hamiltonian PDE. 
It can
be recast as the equation of motion corresponding to a Hamiltonian of a perturbed $\Bbb Z^d$-system of coupled
harmonic oscillators with i.i.d. random frequencies (see (2.2, 2.3)). When $\delta=0$, the linear system has pure
point spectrum: $\sigma(H)=\sigma_{pp}$. This corresponds to the KAM tori scenario. A natural
question is the stability of such invariant tori under small ($0<\delta\ll 1$) perturbations,
which leads to construction of quasi-periodic or almost periodic solutions to (1.1).

\noindent {\it Remark.} Previously in \cite{AF, AFS}, solutions to the nonlinear eigenvalue 
problem
$$(\epsilon\Delta+V)\phi+\delta |\phi|^{2p}\phi=E\phi \qquad \text{on}\quad\ell^2(\Bbb Z^d)$$
were found, which give the time periodic solutions to (1.1) of the particular form 
$$u(j,t)=\phi(j)e^{-iEt}.$$
\smallskip
\noindent {\it A sketch of the construction.}

We expand in the Fourier basis: $e^{in\cdot\omega t}\delta_k(j)$ and as an ansatz, seek
solutions of the form 
$$
u(\ell, t)=\sum_{(j,n)\in\Bbb Z^{d+\nu}}\hat u(j, n) e^{in\cdot\omega t}\delta_j(\ell),\tag 1.8 
$$
with the initial condition 
$$u(\ell, 0)=\sum_{k=1}^{\nu} a_k\delta_k(\ell),\quad\text {satisfying}\, \sum_{k=1}^{\nu}|a_k|\ll 1,\tag
1.9$$ where in (1.9), we identify $\{j_k\}_{k=1}^{\nu}$ with $\{1,...,\nu\}$, $\delta_k$ with
$\delta_{j_k}$ ($k=1,...,\nu$). The unperturbed frequencies are therefore $\omega=\omega(\Cal
V) =\Cal V\in\Bbb R^{\nu}$, where $\Cal V{\overset\text{def }\to =}\{v_{j_k}\}_{k=1}^{\nu}$ are
the random potentials at sites $j_k\in\Bbb Z^d$.

Substituting (1.8) into (1.1), we obtain the following equation for the Fourier coefficients:
$$(n\cdot\omega+\epsilon \Delta_j+V_j)\hat u(j, n)+\delta[(\hat u*\hat v)^{*p}* \hat u](j,n)=0,\tag
1.10$$
where $\hat v(j, n)=\bar{\hat u}(j,-n)$, the convolution $*$ is in the $n$ variable only, $*p$ denotes the
$p$-fold convolution and we added the subscript
$j$ to operators that originated from $\ell^2(\Bbb Z^d)$. We also write the equation for $\hat v$:
$$(-n\cdot \omega+\epsilon\Delta_j+V_j)\hat v(j, n)+\delta[(\hat u*\hat v)^{*p}*\hat v](j,n)=0.\tag
1.11$$
Combining (1.10, 1.11), we then have a closed system of equations for $y=\pmatrix \hat u\\ \hat
v\endpmatrix$, which we write as $$F(y)=0.\tag 1.12$$

Equation (1.12) is a $\Bbb Z^{d+\nu}$ system of equations. Let $y_0=y(t=0)$. 
$$\text{supp }y_0=\{j_k,\,-e_k\}_{k=1}^{\nu}\cup \{j_k,\,e_k\}_{k=1}^{\nu},\tag 1.13$$
where $e_k$ are the unit vectors of $\Bbb Z^\nu$. We seek solutions to (1.12) with
$y$ fixed at the initial condition on $\text {supp }y_0$, i.e., $\hat u(j_k, -e_k)=a_k$, 
$\hat v(j_k, e_k)=\bar a_k$, $k=1,...,\nu$, cf. (2.8).
We make a Lyapunov-Schmidt
decomposition as in \cite{B1,3, CW1,2}. Let $y_0=y(t=0)$. The equations
$$F(y)=0|_{\Bbb Z^{d+\nu}\backslash \text{supp}\,y_0}\qquad \text{on}\quad \ell^2(\Bbb
Z^{d+\nu}\backslash\text{supp}\,y_0)$$
are the so called $P$-equations, the rest are the $Q$-equations. The $P$-equations are used to
determine $y(j, n)$ on $\{\text{supp}\,y_0\}^c$. On supp $y_0$, $y(j, n)$ are held fixed at
the initial condition from (1.9). Instead the $\nu$ $Q$-equations determine $\omega=\omega(\Cal V)$.

We use a Newton scheme to solve the $P$-equations (for more details, see section 3). This leads to
investigate the invertibility of the linearized operators $F_i'(y_i)$, where $y_i$ is the 
$i^{\text {th}}$ approximate solution, $F_i'$ is $F'$ restricted to $[-M^{i+1}, M^{i+1}]^{d+\nu}$
($i\geq 0$) for appropriate $M$.

The random potentials $\Cal V=\{v_{j_k}\}_{k=1}^{\nu}\in\Bbb R^{\nu}$ are the parameters in the
problem. Invertibility of $F_i'(y_i)$ are assured by appropriate incisions in the probability space 
$\Bbb R^\nu$.
Similar to the linear case in \cite{BW}, this is done by using semi-algebraic set techniques
to control the complexity of the sigular sets and Cartan type of theorem for analytic matrix 
valued functions to control the measure.

The main difference with the linear case in \cite{BW} is that $F_i'$ are evaluated at different
$y_i$. But due to rapid convergence of the Newton scheme, made possible by estimates on
$F_{i'}'(y_{i'})$ for $i'<i$, this is within the margin of estimates.

Solving the $P$-equations iteratively is the main part of the work. The solutions to the 
$P$-equations are then substituted into the $Q$-equations to determine $\omega=\omega(\Cal V)$
iteratively by using the implicit function theorem. We obtain time quasi-periodic solutions of the form
(1.8) to (1.1), which are exponentially localized (both in the spatial and Fourier space) to
the initial condition (1.9), with modified frequencies $\omega=\omega(\Cal V)$, which are
($\epsilon+\delta$)-close to the unperturbed frequencies $\Cal V=\{v_{j_k}\}_{k=1}^{\nu}$. 

We therefore have 
\vfill\eject
\noindent {\it Statement of the Theorem.}
\proclaim
{Theorem}
Consider the nonlinear random Schr\"odinger equation
$$i\frac\partial{\partial t}u =(\epsilon\Delta+V)u+\delta |u|^{2p}u,\qquad (p\in\Bbb N^+),\tag 1.14 
$$
where $\Delta$ is the discrete Laplacian defined in (1.2), $V=\{v_j\}_{j\in\Bbb Z^d}$ is a family 
of i.i.d. random variables with common distribution $g$ satisfying (1.3). Fix $j_k\in\Bbb Z^d$,
$k=1,\cdots,\nu$. Let $\Cal R=\{j_k\}_{k=1}^{\nu}\subset \Bbb Z^d$, $\Cal V=\{v_\alpha\}_{\alpha\in\Cal
R}\in\Bbb R^{\nu}$. Consider an unperturbed solution of (1.14) with $\epsilon$, $\delta=0$,
$$u_0(y,t)=\sum_{k=1}^\nu a_ke^{-iv_{j_k}t}\delta_{j_k}(y),$$
with $\sum_{k=1}^\nu |a_k|$ sufficiently small. Let $a=\{a_k\}_{k=1}^\nu$. 

For $0<\epsilon\ll 1$,
$\exists X_\epsilon\subset\Bbb R^{\Bbb Z^d}\backslash\Bbb R^\nu$ of positive probability,
such that for $0<\delta\ll 1$,if we fix $x\in X_\epsilon$, there exists $\Cal G_{\epsilon,\delta}(x; a)\subset\Bbb
R^\nu$, Cantor set of positive measure. There is 
$\omega=\omega_{\epsilon, \delta}(\Cal V;a)$ smooth function defined on 
$\Cal G_{\epsilon,\delta}(x; a)$, such that if $\Cal V\in \Cal G_{\epsilon,\delta}(x; a)$,
then
$$
u_{\epsilon,\delta,x}(y, t)=\sum_{(j,n)\in\Bbb Z^{d+\nu}}\hat u(j, n) e^{in\cdot\omega t}\delta_j(y)
\tag 1.15 
$$
is a solution to (1.14), satisfying 
$$\aligned &\hat u(j_k,-e_k)=a_k,\quad k=1,...,\nu,\\
&\sum_{(j,n)\notin\Cal S}e^{c(|n|+|j|)}|\hat u(j,n)|<\sqrt{\epsilon+\delta},\\
&|\omega-\Cal V|<c(\epsilon+\delta),\endaligned\tag 1.16$$
for some $c>0$, and where $\{e_k\}_{k=1}^\nu$ are the basis vectors for $\Bbb Z^\nu$ and
$\Cal S=\{j_k,-e_k\}_{k=1}^\nu\subset\Bbb Z^{d+\nu}$.
The sets $X_\epsilon$ and $\Cal G_{\epsilon,\delta}(x; a)$ satisfy 
$$\text {Prob } X_\epsilon\to 1,\quad \text{mes } \Bbb R^\nu\backslash \Cal G_{\epsilon,
\delta}(x;a)\to 0\quad
\text{as}\quad \epsilon+\delta\to 0.$$
\endproclaim

\noindent {\it Remark.} The set $X_\epsilon\subset\Bbb R^{\Bbb Z^d}\backslash \Bbb R^\nu$ only depends 
on $\epsilon$; while the set $\Cal G_{\epsilon, \delta}(x; a)\subset\Bbb R^\nu$ depends on $\epsilon$,
$\delta$, $x\in X_\epsilon$ (the random potentials in $X_\epsilon$) and $a$ (the initial amplitude).

\proclaim {Corollary}
There exists $X_{\epsilon, \delta}\subset \Bbb R^{\Bbb Z^d}$ of positive probability,
$0<\epsilon\ll 1$, $0<\delta\ll 1$,
satisfying $$\text {Prob } X_{\epsilon, \delta}\to 1\quad \text{as}\quad\epsilon+\delta\to 0,$$ such
that for initial amplitudes $a$ sufficiently small, there are quasi-periodic solutions to (1.14).
\endproclaim
\smallskip
\noindent{\it Comments on the family of parameters $\{v_j\}_{j\in\Bbb Z^d}$.}

In solving (1.14), we use the basis $e^{in\cdot\omega t}\delta_j$, $(j,n)\in\Bbb Z^{d+\nu}$, (cf.
(1.8)). In the $\Bbb Z^d$ basis $\delta_j$ ($j\in\Bbb Z^d$), the linear operator $H=\epsilon\Delta+V$
is {\it not} diagonalized. Hence $\{v_j\}_{j\in\Bbb Z^d}$ is {\it not} a family of independent 
parameters. This is a slight variation from the ``usual" scenario, where the linear operator 
is diagonalized and the parameters are independent, which is the case in e.g., \cite{B3}.

Here it is convenient to work with the $\Bbb Z^d$ basis $\delta_j$ instead of the basis provided
by the eigenfunctions $\psi_j$ of $H$, as $\psi_j$ depends on $\{v_k\}_{k\in\Bbb Z^d}$. More
precisely, as $\{v_k\}_{k\in\Bbb Z^d\backslash\Cal R}$ is held fixed
on the appropriate probability subspace, $\psi_j$ depends on $\{v_{k'}\}_{k'\in\Cal R}$, which 
serve as parameters for the construction and are therefore varying, (see the statement of
the Theorem).

From the KAM perspective, the normal frequencies are provided by the eigenvalues $\mu_j$ of $H$.
Since $\{v_k\}_{k\in\Bbb Z^d\backslash\Cal R}$ is fixed, the strong localization property (A8)
of $\psi_j$ implies that the normal frequencies $\mu_j$ for $|j|>\rho$, where $\rho$
only depends on the radius of $\Cal R$, can in fact be held fixed. This is close to the
usual terrain, where the normal frequencies are fixed, while the tangential frequencies
vary to avoid small divisors, either via the parameters or via amplitude-frequency modulation
(cf. \cite{B3, KP}).

\noindent{\it An insertion into a larger picture.}

The Theorem presented above is proven for i.i.d. random potentials $V=\{v_j\}_{j\in\Bbb Z^d}$. The
construction used to prove the theorem is, however general. The essential ingredient is a {\it
spectral separation property} on the linearized operator $\tilde H:\,\tilde H=n\cdot\omega+H$,
where $\omega$ are the tangential frequencies, $H$ is the original linear operator (corresponding to
the quadratic part of the Hamiltonian, cf. (2.2)). In the present case, $H=\epsilon\Delta+V$.
Assume $H$ has pure point spectrum and we look at initial conditions localized about the
origin. Below is a tentative formulation of this spectral property.

Let $\mu_j$ be the eigenvalues of $H$. Let $i=(j,n)$, $\lambda_i=n\cdot\omega+\mu_j$ be the
eigenvalues of $\tilde H$. Let $\chi$ be an appropriate function, which depends essentially only on
the initial condition, localized about the origin, $|\chi|\lesssim 1$. Let $\phi_i$, $\phi_i'$
be eigenfunctions of $\tilde H$ (i.e., products of eigenfunctions of $H$ and the exponentials).
Define 
$$K(i, i')=\int\phi_i\chi\phi_{i'}.$$
$\tilde H$ has {\it spectral separation property} if for each scale $L$, $\exists \ell\ll L$, such
that $$|\lambda_i-\lambda_{i'}|\gg K(i,i')\tag 1.17$$
for $\ell\leq|i-i'|\leq L$ ($i\neq i'$).
The $\mu_j$, $\lambda_i$, $\phi_i$ can be replaced by their local version whenever appropriate. 

In the present case, $H=\epsilon\Delta+V$, we use the local version. Assume $\epsilon$ is 
small so that $H$ has A.L. (1.17) is provided by using (A5-7) and restricting to the appropriate probability
subspace, (2.10) and a  direct incision in the frequency space. Related spectral separation properties seem to hold in
\cite{B3, W}. (Compare (1.17) with the nondegeneracy condition in \cite{KP, p164}, where eigenfunctions do not seem
to play an explicit role.)

\noindent{\it Remark.} For the random Schr\"odinger operator $H=\epsilon\Delta+V$ ($\epsilon\ll 1$), no Diophantine
property of the eigenvalues seems to be known at present. So a possible extension of
the standard KAM method, as outlined in e.g., \cite{FSW} is not feasible. It is known following \cite{Mi} however,
that the eigenvalue statistics is Poisson and that in a box of size $N$, the eigenvalue spacing is 
$N^{-p}$ ($p\geq d$). From general considerations, the spectrum of $H$, $\sigma(H)$ is simple
\cite{Si}. 

The construction of time quasi-periodic (or almost-periodic) solutions needs a parameter. This
parameter can sometimes be extracted from amplitude-frequency modulation, see e.g., \cite{KP}.
Nonlinear random Schr\"odinger equation is an equation endowed with a family of parameters, where
the separation property (1.17) can be obtained from A.L. of the linear operator. So it is a
natural candidate for the construction of KAM type solutions. 

The continuum Schr\"odinger equations (linear or nonlinear) are a more frequently studied subject.
The discrete nonlinear Schr\"odinger equation presented here should be seen as the analogue of the
continuum nonlinear Schr\"odinger equation in a compact domain, e.g., on a torus. The $\Bbb Z^d$ lattice
therefore can be seen as the indices of the eigenvalues or eigenfunctions for the underlying linear
Schr\"odinger operator.

Time quasi-periodic solutions have been constructed for the continuum nonlinear Schr\"odinger or
wave equation in 1-D, on a finite interval with either Dirichlet or periodic boundary conditions. See
for example, the works of
Bourgain, Kuksin, P\"oschel and Wayne in \cite{B1, KP, Wy}. In \cite {B3}, time quasi-periodic
solutions are constructed for the 2-D nonlinear Schr\"odinger equation on $\Bbb T^2$. In arbitrary
dimension D, time quasi-periodic solutions for nonlinear Schr\"odinger and wave equations are
treated in \cite{B5, EK}.

The construction presented here is related to those in \cite{B1-5}, which use a Newton scheme
directly on the equations. This direct approach is originated by Craig and Wayne in \cite{CW1,2}. 
It has the advantage of not relying on the underlying Hamiltonian structure. The Hamiltonian
structure does assure, however, the reality of the frequency $\omega$ during the iteration, 
(see section 2 and \cite{B3}).

We end this section by remarking that the present method, as it stands, does not yet extend to
the construction of almost-periodic solutions. This is because the linear equation that serves 
as the starting point of our perturbation is 
$$i\frac{\partial}{\partial t}u=Vu,$$
and {\it not} $$i\frac{\partial}{\partial t}u=(\epsilon\Delta+V)u.$$
In order to construct almost-periodic solutions, we will need more informations on the spectrum
of the linear operator $H=\epsilon\Delta+V$.

In \cite{B2}, the construction of almost-periodic solutions for 1-D nonlinear Schr\"odinger and wave
equations under Dirichlet boundary conditions was made possible by the precise knowledge of the spectrum of the
linear operator and the fact that the perturbation is quartic (in the Hamiltonian). In the present case it is
quadratic. In \cite{B6}, almost-periodic solutions for a 1-D nonlinear Schr\"odinger 
equation under periodic boundary conditions and realistic decay conditions were constructed. In
particular this applies in the real analytic category. Almost-periodic solutions have also been constructed by
P\"oschel \cite{P\"o2} in the case of a nonlinear Schr\"odinger equation, where the nonlinearity is ``nonlocal".

PDE's (such as (1.1)) typically correspond to the so called ``short range" (but not finite range)
case. In the ``finite range" case, which typically corresponds to perturbation of integrable Hamiltonian
systems, almost-periodic solutions have been constructed in e.g., \cite{CP, FSW, P\"o1} among
others.

\bigskip
\head{\bf 2. Hamiltonian representation and Lyapunov-Schmidt decomposition}
\endhead

Recall from section 1, the nonlinear random Schr\"odinger equation
$$
i\frac\partial{\partial t}u =(\epsilon\Delta+V)u+\delta |u|^{2p}u\qquad (p\in\Bbb N^+),\tag 2.1
$$
where $0<\epsilon, \delta\ll 1$, $\Delta$ is the discrete Laplacian as defined in (1.2),
$V=\{v_j\}_{j\in\Bbb Z^d}$ are i.i.d. random variables with common
distribution $g$ as in (1.3). 
The solutions $u=\{u(j,t)\}_{j\in\Bbb Z^d, t\in [0, \infty)}$.

Equation (2.1) can be recast as (infinite dimensional) Hamiltonian equations of motion, with canonical variables
$(u, \bar u)$ and the Hamiltonian
$$
\align
H(u, \bar u)&=\frac{1}{2}\big[\sum_{j, j'\in \Bbb Z^d\times \Bbb Z^d}(\epsilon \Delta+V)_{jj'} u_j\bar
u_{j'} +\bigg(\frac \delta{p+1}\bigg) \sum_j u_j^{p+1}\bar u_j^{p+1}\big]\\
& {\overset\text{def }\to =} H_0(u, \bar u)+\delta H_1 (u, \bar u).
\tag 2.2
\endalign
$$
Equation (2.1) can then be written as
$$
i\dot u=2\frac{\partial H}{\partial\bar u}.\tag 2.3
$$

\noindent
{\it Remark.}
The connection with the usual canonical variables $(p,  q)$ is $u=p+iq, \bar u=p-iq$.
The equation of motion in the $(p, q)$ coordinates is
$$
\dot p=\frac{\partial H}{\partial q}, \dot q =-\frac{\partial H}{\partial p},
$$
which can be rewritten as a single equation (2.3).
(This also explains the factors $i$ and $2$.)

Equations (2.2, 2.3) show that (2.1) can be viewed as a perturbed $\Bbb Z^d$-system of coupled harmonic
oscillators with i.i.d. random frequencies.
The perturbation $H_1$ can be of a more general type, e.g.,
$$
H_1(u, \bar u)=\sum_{j, j '\in \Bbb Z^d\times \Bbb Z^d} a_{jj'}u_j^{p+1}\bar u_{j'}^{p+1}\tag 2.4
$$
with $a_{jj'} = a_{j'j} $ decaying exponentially or polynomially of sufficiently high degree as $|j-j'|\to\infty$.
The reason we mention (2.4) is to stress that the construction we present below {\it does not} rely on 
integrability of the system.
It also carries through for $H_1$ of type (2.4), although we only present it for $a_{jj'}=\delta_{jj'}$.

The goal of the rest of the paper is to seek time quasi-periodic solutions to (2.1) for appropriately
chosen {\it localized} initial conditions.
We hence expand $u$ in the  basis
$$
e^{in\cdot \omega t}\delta_k(j),\tag 2.5
$$
where $n\in \Bbb Z^\nu, \omega\in\Bbb R^\nu$, $k, j\in\Bbb Z^d$, $\delta_k(j)$ is the canonical basis for $\Bbb Z^d$.
$\delta_k(j)$ is a natural basis here due to smallness of $\epsilon$.
(In [B1-3], the spatial basis is given by the eigenfunctions of the linear operator.
The $k$-labeling there is the eigenvalue labeling.)

In the basis (2.5), (2.1) becomes
$$
(n\cdot \omega+\epsilon \Delta_j+V_j)\hat u(j, n)+\delta\widehat{\frac{\partial H_1}{\partial \bar u}}
(j, n)=0,\tag 2.6
$$
where $n\in \Bbb Z^\nu, j\in\Bbb Z^d$, $H_1$ as defined in (2.2) and $\hat u$ are the Fourier
coefficients of $u$:
$$
u(k, t)=\sum_{(j,n)}\hat u(j, n) e^{in\cdot\omega t}\delta_j(k).\tag 2.7
$$
We have also put the subscript $j$ on operators that operate in the spatial $(\Bbb Z^d)$ variable only.
(This is the same notation as in [BW].)

In view of the Theorem, we seek solutions to (2.6) with the constraint
$$
\hat u(j_k, -e_k)=a_k \qquad (k=1, \ldots\nu),\tag 2.8
$$
where $j_k\in \Bbb Z^d$, $e_k$ are unit vectors in $\Bbb Z^\nu$, $a_k$ are fixed. 
Assume $\omega_1,\,\omega_2,...,\,\omega_\nu$ are rationally independent, i.e., $\omega=\{\omega_i\}_{i=1}^\nu\in\Bbb
R^\nu$ is a Diophantine vector, which will be the case when the Theorem applies. A time shift and a limiting
argument (since the Kronecker flow is dense) permit us to assume $a_k$ are real. Hence from now on $a_k\in\Bbb R$,
$k=1,...,\,\nu$.

Due to the smallness of
$\epsilon$, we take our initial unperturbed linear equation to be
$$
i\frac{\partial}{\partial t}u=Vu.\tag 2.9
$$
The conditions in (2.8) thus correspond to the initial unperturbed solution
$$
u_0(k, t) =\sum^\nu_{\ell=1} a_\ell e^{-iv_{j_\ell} t}\delta_{j_\ell}(k).\tag 2.10
$$

Let 
$$
a=\{a_k\}^\nu_{k=1}\in\Bbb R^\nu,\quad \Cal R =\{j_k\}^\nu_{k=1} \subset \Bbb Z^d, \quad \Cal V
=\{v_\alpha\}_{\alpha\in\Cal R}\in\Bbb R^\nu.$$ We constructively show that for $\epsilon$ small
enough, $\exists X_\epsilon\subset \Bbb R^{\Bbb Z^d}\backslash
\Bbb R^\nu$ of positive probability, satisfying $\text {Prob}\,X_\epsilon\to 1$, as $\epsilon\to 0$,
such that if we fix $x\in X_\epsilon$, for $\delta$, $a$ small enough,
there exists
$\Cal G_{\epsilon, \delta} (x; a)\subset \Bbb R^\nu$, Cantor set of positive measure, satisfying
$\text {mes}\,\Bbb R^\nu\backslash\Cal G_{\epsilon, \delta} (x; a)\to 0$, as $\epsilon+\delta\to 0$. We can find
$\omega=\omega (\Cal V; a)$, smooth function defined on $\Cal G_{\epsilon, \delta} (x; a)$ and $\hat
u$ such that (2.6) holds.
$\omega$ and $\hat u$ are determined simultaneously in an iterative way.

Toward that end, we first perform a Lyapunov-Schmidt type decomposition (see [B1-3, CW1,2]) of (2.6).
Let
$$
\Cal S=\{(j_k, -e_k)|k=1, \ldots\nu\} \subset \Bbb Z^{d+\nu}.\tag 2.11
$$
From (2.10), $\Cal S=\text {supp } u_0$, $u_0$ is a solution to (1.14) when $\epsilon=\delta=0$. We
call
$\Cal S$ the resonant set and consider the $\nu$ equations
$$
\aligned
[(n\cdot\omega+\epsilon\Delta_j+V_j)\hat u](j_k, -e_k)+\delta\, \widehat{\frac{\partial H_1}{\partial\bar
u}}\quad &(j_{k},-e_k)=0\\
 &(k=1, \ldots, \nu)
\endaligned
\tag 2.12
$$
obtained by taking $(j, n) \in \Cal S$.

They form the finite system of $Q$-equations.
The remaining infinite system of equations are called the $P$-equations
$$
(n\cdot \omega+\epsilon \Delta_j+V_j) \hat u(j, n)+\delta \ \widehat{\frac{\partial H_1}{\partial \bar
u}} (j, n) =0,\quad (j, n)\not\in \Cal S.\tag 2.13
$$
The $P$-equations are used to determine $\hat u(j, n)$ for $(j, n)\not\in \Cal S$.
(Recall from (2.8) that $\{\hat u (j, n), (j, n)\in \Cal S\}=a$ are given.)

Once $\hat u(j, n)$ are determined, the $Q$-equations in (2.12) are used to determine $\omega=\omega(\Cal V, a)$ via
the implicit function theorem.
Since $a$ is real, $H_1$ is a polynomial in $u, \bar u$ with real coefficients, the solution $\hat u$ to (2.13) will
be real and hence also $\omega$ determined from (2.12).
(For more details, see the comment after (2.3) of [B3].)

To solve (2.13), we duplicate the equation for $\bar u$ to form a closed system.
Let 
$$
\align
v&=\bar u\\
\hat v(j, n)&=\bar{\hat u}(j, -n)\\
-\Cal S&= \{(j_k, +e_k)|k=1, \ldots\nu\}\subset \Bbb Z^{d+\nu}.\tag 2.14
\endalign$$
(The flip in  sign in the second equation of (2.14) is solely in order that the
convolution coming from the nonlinearity obeys the usual sign convention.)

We then have the closed system of $P$-equations
$$
\cases (n\cdot \omega+\epsilon \Delta_j+V_j)\hat u(j, n)+\delta\ \widehat{\frac{\partial H_1}{\partial v}}\, (j,
n)=0,\quad (j, n)\not\in \Cal S,\\
(-n\cdot\omega+\epsilon\Delta_j+V_j)\hat v(j, n)+\delta \ \widehat{\frac{\partial H_1}{\partial u}} \,
(j, n)=0,\quad (j, n)\not\in -\Cal S.
\endcases
\tag 2.15
$$
For $H_1$ as in (2.2), (2.15) takes the explicit form
$$
\cases
[(n\cdot\omega+\epsilon \Delta_j+V_j)\hat u](j, n)+\delta[(\hat u*\hat v)^{*p}* \hat u](j,n)=0,\\
[(-n\cdot \omega+\epsilon\Delta_j+V_j)\hat v](j, n)+\delta[(\hat u*\hat v)^{*p}*\hat v](j,n)=0,
\endcases \quad (j, n)\in \Bbb Z^{d+\nu} \backslash (\Cal S\cup -\Cal S)\tag 2.16
$$
where the convolution $*$ is in the $n$ variable only.
We solve (2.16) by using a Newton iteration scheme to be amplified in the next section.
We also identify $\hat u$ with $u$, $\bar{\hat u}$ with $v$ and write $y$ for $\pmatrix \hat u\\ \hat
v\endpmatrix$.

\head
{\bf 3. Newton scheme}
\endhead

Let $F$ denote the left hand side (LHS) of (2.16).
Our task is to restrict the set of $(\omega,\Cal V)\in\Bbb R^{2\nu}$ in order to find $y$ such that
$$
F(y)=0,\tag 3.1
$$
so that (2.16) is resolved. We use a Newton iteration. Recall first the formal scheme.

Starting from the initial approximant $y_0$, solution to (1.14) and its conjugate when
$\epsilon=\delta=0$, the successive approximants
$y_i$ are defined by
$$
\Delta_{i+1} y  \ \ {\overset \text{\rm def }\to =} \ \ y_{i+1} -y_i=-[F'(y_i)]^{-1}F(y_i).\tag 3.2
$$
Let $T$ denote the linearized operator $F'$.
From (2.16)
$$
T=D+\delta S,\tag 3.3
$$
where $D$ is diagonal (in the $n\in \Bbb Z^\nu$ variables)
$$
\align
D &=\pmatrix n\cdot\omega+\epsilon \Delta_j+V_j&0\\ 0& -n\cdot\omega+\epsilon \Delta_j+V_j\endpmatrix\\
&  \ {\overset\text{def }\to=} \ \pmatrix D_+&0\\ 0&D_-\endpmatrix\tag 3.4
\endalign
$$
and
$$
S=S(u, v)=\pmatrix (p+1)(u*v)^{*p}& p(u*v)^{*p-1}*u*u\\ p(u*v)^{*p-1}*v*v & (p+1)(u*v)^{*p}\endpmatrix
\quad (p\in\Bbb N^+)\tag 3.5
$$
evaluated along the previous approximant.
We note that $S$ is self-adjoint, although this does not play a role in our construction.

Denote by $\Vert \ \Vert$ the $\ell^2$ norm of a vector or operator on $\ell^2 (\Bbb Z^{d+\nu})$.
Using (3.2), the error of the approximation at stage $(i+1)$ can be estimated
$$
\align
F(y_{i+1})&= F(y_i)+F'(y_i)(y_{i+1}-y_i)\\
&\qquad +\Cal O (\Vert y_{i+1}-y_i\Vert^2)
\\
&=\Cal O (\Vert y_{i+1}-y_i\Vert^2).
\tag 3.6
\endalign
$$
So using (3.2)
$$
\Vert F(y_{i+1})\Vert =\Cal O (\Vert[F' (y_i)]^{-1}\Vert^2)\Vert F(y_i)\Vert^2.\tag 3.7
$$
The crux of the matter is thus to control $\Vert[F'(y_i)]^{-1}\Vert$ in order that
$$
\Vert F(y_{i+1})\Vert\ll \Vert F(y_i)\Vert.\tag  3.8
$$
(Note the squaring of the norm of $F(y_i)$ in the RHS of (3.7), which makes this feasible.)

Since (3.1) represents an infinite system of equations and the initial condition (2.10) is localized in a compact
region in $\Bbb Z^{d+\nu}$, to control the norm of $[F' (y_i)]^{-1}$ we implement the Newton scheme in a slightly
modified way, gradually increasing the size of the system that we consider.

Let $M\in \Bbb N^+$, which can be assumed large in order that $[-M,M]^{d+\nu}\supset 2p\text{ supp\,}
y_0$, in view of (3.4, 3.5) (See also (3.10) below.).
At stage $i$, let $N=M^{i+1}$ and $T_N(y_i)$, the restriction of  $T(y_i)$ to $[-N, N]^{d+\nu}$.
We define
$$
\Delta_{i+1} y=y_{i+1} -y_i \ {\overset\text{def }\to=} \ -\big(T_N(y_i)\big)^{-1} F(y_i).\tag 3.9
$$
So
$$
\align
F(y_{i+1})&=F(y_i)+F'(y_i)(y_{i+1}-y_i)+\Cal O(\Vert y_{i+1}-y_i\Vert^2)\\
&=(T-T_N)(y_{i+1}- y_i)+\Cal O(\Vert y_{i+1}-y_i\Vert^2)\\
&=-[(T-T_N)T_N^{-1}]F(y_i)+\Cal O(\Vert T_N^{-1}\Vert^2\Vert F(y_i)\Vert^2),\tag 3.10
\endalign
$$
where we used (3.9).
Compared to (3.6) the first term in the RHS of (3.10) is new. Moreover it is only {\it linear} in 
$F(y_i)$. 
This necessitates the control of off-diagonal decay of $T$ and $T_N^{-1}$ evaluated at $y_i$, in
addition to that of $\Vert T_N^{-1}\Vert$.

\noindent {\it  The control of} $T_N^{-1}$

Recall that the $0^{\text{th}}$ approximant (initial condition) to (3.1) $y_0$ has compact support, 
$\text{supp\,} 
y_0=\Cal S\cup -\Cal S$, where $\Cal S$ and $-\Cal S$ are defined in (2.11, 2.14).
From (3.3-3.5), $T(y_0)$ is a diagonal dominated matrix with finite range off-diagonal elements.
So off-diagonal decay of $T(y_0)$ is automatically satisfied.

Assume the successive approximants $y_i$ are {\it uniformly} (in $i$) exponentially localized about 
$\Cal S\cup-\Cal S$ (cf. (1.16)).
This assumption will be verified later from the construction itself in view of (3.9, 3.10).
From (3.5) the successive $S(y_i)$ have {\it uniformly} exponentially decaying off-diagonal elements
in the $n$ direction, and are diagonal in the $j$ direction, with a prefactor which decays exponentially
in
$j$. (The exponential decay of the prefactor stems from the uniform exponential decay assumption on
$y_i$.) Hence
$T(y_i)$ are of the type (although more complicated) of matrix operator studied in [BW].

To study the $T$'s, we introduce, as in [BW] an auxiliary parameter $\theta \in \Bbb R$.
We consider instead
$$
T^\theta =D^\theta+\delta S, \tag 3.11
$$
where
$$
D^\theta =\pmatrix n\cdot\omega+\theta+\epsilon\Delta_j+V_j&0\\
0& -(n\cdot\omega+\theta)+\epsilon \Delta_j+V_j\endpmatrix
$$
$$
{\overset\text{def }\to=} \ \pmatrix D_+^\theta&0\\
0&D^\theta_-\endpmatrix\tag 3.12
$$
and $S$ is as before in (3.5).
Similarly we define $T_N^\theta(y_i)$, where $N=M^{i+1}$ as in (3.9).

In section 4, fix $x$ in a good set of probability space, where there is Anderson localization for the
linear random Schr\"odinger operator $H_j=\epsilon\Delta_j+V_j$, so that (A6) holds. (For precise
details see the appendix.) Assuming $\omega$ Diophantine, $y_i$ uniformly (in $i$) exponentially
localized about $\Cal S\cup -\Cal S$, we bound the norm of $[T_N^\theta(y_{i'})]^{-1}$, where
$N=M^{i}$, $i'>i$ as in (4.9) (The precise relation between $i$ and $i'$ is dictated by the
construction in section 5.) and establish exponential decay properties of its off-diagonal matrix
elements on a set of $\theta$ of small complementary measure. 

In Lemma 4.1, fix any $y_k$, we bound $[T_N^\theta(y_k)]^{-1}$ for all $N$. We then use it to obtain estimates on 
$[T_N^\theta(y_{i'})]^{-1}$, where $N=M^i$, $i'>i$ satisfying the restriction in the third line of (4.9).
This bound is 
abstract in the sense that $\omega$, $\Cal V$, $y_{i'}$ are viewed as {\it independent}
parameters for the time being.

As in [BW], this is an iteration process, using semi-algebraic set techniques and
Cartan-type theorem. To start the iteration, we neglect $\delta S$ and exclude a set of $\theta$
such that $D_N^\theta=D_M^\theta$ has a small diagonal element. To continue the iteration,
we also need to exclude a set in $\omega$ of small measure. It is important to remark that this set in
$\omega$ is {\it independent} of $\Cal V$, $y_{k}$. It only depends on $x\in\Bbb R^{\Bbb
Z^d}\backslash\Bbb R^\nu$. We stress that for fixed $x$ in the good probability set, Lemma
4.1 holds for {\it any} fixed $\omega$ in the good frequency set, {\it any} $\Cal V\in\Bbb R^{\nu}$ and {\it any} $y_{k}$
which satisfy (H1-3) in section 4. 
The set of excluded $\theta$, $\Cal B(N)$ depends, of course, on $x$, $\omega$, $\Cal V$,
$y_{k}$. 

In section 5, we iteratively transfer the estimates on $T^\theta_N(y_{i'})$ in $\theta$ into estimates on 
$T_{\bar N}(y_{\bar i})$ in
$(\omega, \Cal V)$, where $N=M^{i+1}$, $\bar N=M^{\bar i+1}$, $\bar i>i'>i$ to be made precise, using
the resolvent equation and taking into account
the $Q$-equations, which are implicit functions relating 
$\omega$, $\Cal V$, $y_i$.
(Recall that $\theta$ is an auxiliary variable.
In the original problem (3.3), $\theta$ is fixed at 0.) $x$ is fixed in the good set of probability
space as in section 4.

For the first $i$ iterations, we treat $\delta S$ as perturbation and use a direct $\epsilon$, $\delta$
perturbation series.
Instead of excluding a set of $\theta$ as in section 4, we exclude a set of $(\omega, \Cal V)\in \Bbb R^{2\nu}$,
so that in the complement, $T_N(y_i)$ are invertible with exponentially decaying off-diagonal elements.

This generates an initial set of ``good" intervals: $\Bbb
R^{2\nu}\supset\Lambda_1\supset\Lambda_2\supset\cdots\supset\Lambda_i$ in the ($\omega$, $\Cal V$)-space.
Using the Newton scheme in (3.9), this also shows that $y_0$, $y_1$, $\cdots$, $y_i$ are exponentially
localized about $\Cal S\cup -\Cal S$. (Recall that $y_0$ is the initial condition, $\text {supp}\,
y_0=\Cal S\cup -\Cal S$.)

Starting from the ($i+1$)$^{th}$ iteration, aside from direct $\delta$ perturbation series,
for certain parts of the estimates, (which concerns the regions far from the origin in the $\Bbb Z^d$
direction), we need to keep 
$\delta S$. This is the heart of the matter. In technical terms, we need to deal with more general semi-algebraic
sets, which are not solely defined by products of monomials. For such semi-algebraic sets,
we use $Q$-equations and a decomposition lemma (Lemma 9.9 in \cite {B5}, restated as Lemma 5.3 in
section 5) to transfer the measure estimates in $\theta$ in Lemma 4.1 into measure estimates in
$\omega=\omega(\Cal V)\in\Bbb R^{\nu}$. Using perturbation, this in turn generates a new set of
intervals $\Lambda_{i+1}\subset
\Lambda_{i}\subset\cdots\subset\Bbb R^{2\nu}$, in the $(\omega,\Cal V)$ space, on which $T_{M^{i+1}}(y_i)$
is invertible and whose inverse has uniformly (in $i$) exponentially decaying off-diagonal elements.

In section 6, using the Newton scheme (3.9), we construct $y_{i+1}$. The (uniformly in $i$) exponential
localization about $\Cal S\cup -\Cal S$ is preserved. Hence Lemma 4.1 is now available at $y_{k}=y_{i+1}$
for future iterations. We generate $\Lambda_{i+2}$, $y_{i+2}$...

Section 7 summarizes the entire construction. It is merely meant as a recapitulation of the sequence
of events leading to the theorem.

\noindent {\it Two technical subtleties}

\itemitem{$\bullet$}
The estimates in section 4 are obtained along the construction devised in [BW].
However for the application to (4.9), 
$T^\theta_N(y_{i'})$ need to be evaluated at different $y_{i'}$ at different scales $N$.
Due to the {\it uniform} exponential decay estimates on $y_{i'}$, Lemma 4.1 can be applied as explained after
the statement of Lemma 4.1.

\itemitem{$\bullet$} From the $P$-equations, the $y_i$'s are constructed on a {\it good} set of $(\omega, \Cal V)\in
\Bbb R^{2\nu}$.
(This set eventually becomes a Cantor set.)
On the same set of $(\omega, \Cal V)$ we also have estimates on $\partial y_i$, where $\partial$ is with respect to
$\omega$ or $\Cal V$.
Using this, we can construct $y_i$ which is smoothly defined on the whole parameter space $(\omega, \Cal V)$.
(Note that outside the good set in ($\omega$, $\Cal V$), $y_i$ is no longer close to a solution to $F(y_i)=0$.)
Substituting for $y_i$, the $Q$-equations are therefore defined on the whole parameter space $(\omega, \Cal V)$.
We can then use the standard implicit function theorem to determine $\omega=\omega(\Cal V)$.

\head
{\bf 4. $P$-equations and statement in $\theta$}
\endhead

Recall the system of $P$-equations in (2.16)
$$
\cases
(n\cdot \omega + \epsilon\Delta_j+V_j)\hat u +\delta (\hat u*\hat v)^{*p}*\hat u=0\\
(-n\cdot \omega+\epsilon\Delta_j+V_j)\hat v+\delta (\hat u*\hat v)^{*p}*\hat v=0
\endcases
\tag 4.1
$$
on $\ell^2\big(\Bbb Z^{d+\nu}\backslash (\Cal S \cup -\Cal S)\big)$,  where $\Cal S$, $-\Cal S$ are as defined in (2.11,
2.14), which are collectively written as $F(y)=0$, with $y=\pmatrix\hat u\\ \hat v\endpmatrix =\pmatrix u\\
v\endpmatrix$.

We solve (4.1) using a Newton scheme, with the family of linearized operators $T(y_i)$, evaluated at the
$i^{\text{th}}$ approximant $y_i$
$$
T(y_i)=D+\delta S(y_i)\tag 4.2
$$
where
$$
\align
D&=\pmatrix 
n\cdot\omega+\epsilon \Delta_j+V_j&0\\ 0& -n\cdot \omega+\epsilon\Delta_j+V_j\endpmatrix\tag 4.3
\\
&=\pmatrix D_+ &0\\ 0&D_-\endpmatrix
\endalign
$$
and
$$
S(y_i)=\pmatrix (p+1)(u_i*v_i)^{*p}& p(u_i*v_i)^{*p-1}*u_i*u_i\\
p(u_i*v_i)^{*p-1} *v_i*v_i&(p+1)(u_i*v_i)^{*p}\endpmatrix
\tag 4.4
$$
as in (3.4, 3.5).

In view of the Newton scheme in (3.9), we need to study the family of restricted operators: 
$T_N(y_i)$, $N=M^{i+1}$, $M$ assumed large depending on $p$,
$$
T_N(y_i)= R_NT(y_i)R_N\tag 4.5
$$
and $R_N$ is the characteristic function of the set $[-N,N]^{d+\nu}$. This will be achieved in section 5 by using
the resolvent identity, covering $[-M^{i+1},M^{i+1}]^{d+\nu}$ with the interval 

\noindent $I=[-M^{i},M^{i}]^{d+\nu}$
and smaller intervals $J=[-M_0,M_0]^{d+\nu}+k$, $\frac{1}{2}M^i<|k|<M^{i+1}$, $M_0\sim(\log N)^{C/2}$, $\qquad$
(4.6)

\noindent (see (5.5)), and restricting the set of $(\omega, \Cal V)\in\Bbb R^{2\nu}$.

Toward that end, as previously mentioned in section 3, we introduce an additional parameter $\theta\in\Bbb R$ and
let
$$
D^\theta=\pmatrix D_++\theta&0\\ 0&D_--\theta\endpmatrix,
$$
$$
T^\theta(y_i)= D^\theta+\delta S(y_i),
$$
$$T_N^\theta(y_i)= R_NT^\theta(y_i)R_N\tag 4.7$$
As mentioned there, we temporarily view $\omega$, $\Cal V\in\Bbb R^{\nu}$ as independent parameters
in this section. In the same vein, we also dissociate $y_i$ from $\omega$, $\Cal V$, assuming only that 
they satisfy (H1-3) below. It is only in section 5 that we restrict to $\omega=\omega(\Cal V)$, determined
from the $Q$-equations and $y_i$, the $i^{th}$ approximate solution to (4.1), which depends on $\omega$, $\Cal V$.

In the rest of the section
$\omega, \Cal V$ are held {\it fixed}, only $\theta$ is varying.
Note that $T^{\theta}$ is of the form $T^{\theta}\ {\overset\text{def }\to=}
T(\theta)=T'(\theta+n\cdot\omega)$. Later in section 5, we transfer the estimate in $\theta$ into
estimates in $\omega$, hence $\Cal V$ by restricting $\theta$ to be of the form $\theta=n\cdot\omega$,
thereby resolving (4.1), (which is at $\theta=0$) on the good set of $\omega, \Cal V$.

Newton scheme is an iterative scheme. The estimate on $[T_I(y_i)]^{-1}$ for $I$ defined in (4.6) is easily obtained by
perturbation arguments on  $[T_I(y_{i-1})]^{-1}$ known from the previous step, which is the step to construct
$y_{i}$, see (3.9). The main task is to estimate
$[T_J(y_i)]^{-1}$ for $J$ defined in (4.6). We therefore study $[T^\theta_{[-M_0,M_0]^{d+\nu}}(y_i)]^{-1}$ (and
later in section 5, we restrict to $\theta=k\cdot\omega$, $\frac{1}{2}M^i<|k|<M^{i+1}$).  This is the subject of 
Lemma 4.1and its application. Note that $M_0$ corresponds to the size of interval at a stage $[\log M_0/\log M]$, which procedes
$i$, while the linearized operator $T$ is evaluated at $y_i$: $T=T(y_i)$.

Assume $y_i$ satisfies
\itemitem{(H1)}
$ \text{supp } y_i\subseteq [-M^i, M^i]^{d+\nu}\qquad (i\geq 1).
$

(This is by construction, see (3.9).)

\itemitem{(H2')} \ $\Vert \Delta_iy\Vert =\Vert y_i-y_{i-1}\Vert <\sqrt{\epsilon+\delta} M^{-b^i}\qquad
 (i\geq 1)$

for some $1<b<2$ in view of (3.6, 3.9) and $b$ will be specified in (6.20). (Recall $y_0=\pmatrix 
u_0\\  v_0\endpmatrix$, $u_0$ defined in (2.10), $v_0=\bar u_0$.)
 \itemitem{(H3)}  \ $|y_i(k)|< e^{-\alpha|k|} \quad \text { for some $\alpha>0$ (uniform in $i$)}.
$

(There is no constant in front of the exponential, as we assume small initial data:
$|a_\ell|<1, \ell =1, \ldots\nu$. See (1.9).)

\noindent{\it Remark.} Using (3.6) in (3.9)
$\Vert\Delta_iy\Vert^2<\Vert\Delta_{i+1}y\Vert<\Vert\Delta_iy\Vert$, assuming appropriate
condition on $T_N^{-1}$. This is consistent with $1<b<2$ in (H2').

\itemitem{(H1-3)} will be verified along the iteration later in sections 5, 6 using Lemma 4.1 below.

Let $\Lambda_N =[-N, N]^d, X_N\subset \Bbb R^{\Lambda_N}\backslash \Bbb R^{\nu}$
be a set, where $\epsilon \Delta_j+V_j$ has A.L. at scale $N$, in a sense to be made
explicit in the process of the proof,
$X_N$ is asymptotically $(\epsilon\to 0)$ of full measure. (Recall also that $\Cal
V=\{v_{j_k}\}^\nu_{k=1}\in
\Bbb R^\nu$ with measure $\prod^\nu_{k=1} g(v_{j_k}) dv_{j_k}$, is the parameter set.)

\noindent {\it Definition}

For $A$, $c>0$,  $DC_{A,c}(N)
\subset \Bbb R^\nu$ is the set such that if $\omega\in DC_{A, c}(N)$, then
$$
\Vert n\cdot\omega\Vert_{\Bbb T} \geq \frac{c}{|n|^A},  
\ n\in [-N, N]^\nu\backslash\{0\}.\tag 4.8
$$
$\quad DC_{A, c}\subset \Bbb R^\nu$ is the set of $\omega$ such that (4.8) is satisfied for all $N$.

$\Cal B_{\beta,\gamma}(N)\subset\Bbb R$ is the complement of the set of $\theta\in\Bbb R$ defined by  
$$
\aligned
&\Vert T_N^\theta(y_{i'})^{-1} \Vert < e^{N^\beta}\qquad (0<\beta<1),\\
&|[T_N^\theta(y_{i'})]^{-1}(k, k')|<e^{-\gamma|k-k'|}\qquad (\gamma>0),\\
&\text{for }|k-k'|>N/10, N=M^{i+1}, 
i'\simeq \frac{\log M}{\log b} i>i \text{ ($b$ as in (H2'))},\endaligned\tag 4.9
$$
$y_{i'}$ satisfies (H1-3); $i'$ is chosen in view of later construction in Lemma 5.1, (see in particular (5.14,
5.15)). This means of course that at least the first $\Cal O(\frac{\log M}{\log b})$ approximants are obtained by
direct perturbation series in $\epsilon, \delta$. So $\alpha=\Cal O(1)|\log (\epsilon+\delta)|$. In general we write $\Cal B(N)$ for  $\Cal B_{\beta,\gamma}(N)$, unless the parameters $\beta$, $\gamma$ need to be emphasized.

\noindent{\it Remark.} At this stage of the construction, it is sufficient to have a lower bound on $b$. This can
be easily obtained in the first few perturbation series by adjusting $\epsilon$, $\delta$. For later purposes, we 
mention that the Diophantine condition (4.8) will be used for $\omega=\omega_{i'}$,   
the ${i'}^{th}$ approximation.

The inequalities in (4.9) are proven iteratively as in \cite{BW}. The rate of decay $\gamma$ will
deteriorate with iteration. So $\gamma=\gamma_N$. But the decrease will decrease with increasing scales and we have 
$\gamma_N>\gamma/2$ for all $N$, cf. Lemma 4.1 and the paragraph after. This rate of decay determines 
the rate of decay of $y_i$. So this is consistent with the assumption (H3).

Inspecting the definition of $D_\theta$ in (4.7, 4.3),
$\theta$ is {\it not} equivalent to a spectral parameter. Hence we need to resort to Cartan-type
theorems as in the wave case in \cite{BW}. This necessitates that we obtain estimate (4.9) for more
general regions than cubes at each scale $N$, the elementary regions to be defined below.

\noindent{\it Remark.} The various approximants $y_i$ are still evaluated using cubes $[-N,N]^{d+\nu}$,
$N=M^{i+1}$ as in (3.9). It is only at each scale $N$ that we also look at restrictions of $T(y_{i'})$
($i'\simeq \frac{\log M}{\log b} i$) to these more general regions.
\smallskip
\noindent{\it Elementary Regions}

An {\it elementary region} is a set $\Lambda$ of the form
$$
\Lambda {\overset\text{def }\to=} R\backslash (R+k), \, k\in\Bbb Z^{d+\nu} \text { is arbitrary}\tag 4.10
$$
and $R$ is a hyper-rectangle
$$
\align
R=\{\ell'\in \Bbb Z^{d+\nu}| |&\ell_i'-\ell_i|\leq N_i, i=1, \ldots d,\,d+1,..., d+\nu,\\
&\ell' =\{\ell_i'\}_{i=1}^{\nu+d} \in \Bbb Z^{d+\nu},
 \ell =\{\ell_i\}_{i=1}^{\nu+d} \in\Bbb Z^{d+\nu}\}
\tag 4.11
\endalign
$$
Let $N=\max_i N_i {\overset\text{def }\to=} N_{i_{\max}}$. 
Assume $\ell \in \Bbb Z^{d+\nu}$ fixed.
We call $\ell$ the center of $R$. $\Cal E \Cal R(N)$ (at a fixed center) is then defined to be the set of all
regions obtained by varying
$k\in \Bbb Z^{d+\nu}$ and $N_i\,(i\not= i_{\max})$, keeping $N_i\leq N$.
We say $2N$ is the diameter of the elementary regions.

To be economical, we extend the notation $T_N$ to mean $T_{\Lambda(N)}=R_{\Lambda(N)}TR_{\Lambda(N)}$, for any 
$\Lambda(N)\in\Cal E\Cal R(N)$ and where $R_{\Lambda(N)}$ is the characteristic function of the set $\Lambda(N)$; 
$\Cal B(N)$ is then the corresponding bad set, on which (4.9) are violated. For the purpose of constructing
approximate solutions, we only need to specialize to $N=M^{i+1}$, ($i\geq 0$). However to state the various
intermediate technical lemmas, it is more convenient to let $N$ be any integer. 

Fix any $y_k$ satisfying (H1-3). Let $T_N$ be the linearized operator evaluated at $y_k$ for all $N$,
i.e., $T_N=T_N(y_k)$. With  a slight abuse of notation, we also let $\Cal B_{\beta,\gamma}(N)$ be the 
corresponding bad set. 
The main goal of this section is to prove
\proclaim
{Lemma 4.1} Fix $0<\sigma<\frac{1}{6(d+\nu)}$, $\sigma<\beta<1$, $\bar N_0$ sufficiently
large, $\max(\frac{1}{\sigma}, 6(d+\nu))\ll C<\bar N_0^{\sigma/2}$.
There exist $\epsilon_0>0$, $\delta_0>0$, such that for all $0<\epsilon<\epsilon_0$, $0<\delta<\delta_0$, there
exists $X\subset\Bbb R^{\Bbb Z^d}\backslash\Bbb R^\nu$ with 
$$\text {\rm mes\,}X\geq 1-\Cal O(1)\bar N_0^{-\kappa},\tag 4.12$$
where $\kappa=\kappa(C, p', d)>0$, $p'$ as in (A2).
Fix $x\in X$, $\exists \Omega\subset\Bbb R^{\nu}$ (independent of $\Cal V\in\Bbb R^\nu$, $y_{k}$), with
$$\text {\rm mes\,}\Omega\leq e^{-\bar N_0^{\kappa'}},$$
where $\kappa'=\kappa'(C,\beta)>0$, such that if 
$$
\omega\in DC_{A, c}\backslash \Omega\tag 4.13
$$
then for any $\Lambda(N)\in\Cal E\Cal R(N)$, $N\geq\bar N_0$
$$
\text{\rm mes\,} \Cal B_{\beta,\gamma_N}(N)\leq e^{-{N}^{\sigma}}.\tag 4.14
$$
where $\gamma_N\geq \alpha-\bar N_0^{-\delta'}$, ($\delta'>0$), for all $N$, $\alpha$ as in (H3), $\alpha=\Cal O(1)|\log(\epsilon+\delta)|$.
\endproclaim
\noindent
{\it Remark.}
$\Cal B_{\beta,\gamma_N}(N)$ depends only on $y_{k}$, $\omega, \Cal V$ as $x$ is fixed. 
In the proof of Lemma 4.1, only (H3) on $y_k$ is needed.

In order to obtain (4.9) at all scales, we apply Lemma 4.1 as follows. From the third line of (4.9), for any 
fixed $y_k$, we only need the lemma at scale $N=M^i$ with $i=\frac{\log b}{\log M}k$. To go to scale 
$N'=N^C=M^{iC}$ with the corresponding $y_{k'}$, $k'=\frac{\log M}{\log b}iC$, we first use (H2'), which 
gives that   
$$
\Vert T_N^\theta(y_{k})-T_N^\theta(y_{k'})\Vert \leq\Cal O_{d+\nu}(1)\Vert y_{k} -y_{k'}\Vert
\leq \Cal O_{d+\nu}(1)M^{-b^{k}}.\tag 4.15
$$
(4.15) shows that we have essentially the same estimate on $\Vert T_N^\theta(y_{k'})\Vert$ as for  $\Vert T_N^\theta(y_{k})\Vert $.

We use $N$ as the initial scale instead of $\bar N_0$, the proof of the inductive step in 
Lemma 4.1 then gives instead
$$\gamma_{N'}(y_{k'})=\gamma_{N}(y_{k})-N^{-\delta'}, (\delta'>0),\, N'=N^C.\tag 4.16$$
From (3.9), the decay rate $\alpha_{i+1}$ of $y_{i+1}$ is governed by the decay rate of
$[T_N(y_i)]^{-1}$, where $N=M^{i+1}$. (Note that $\theta$ is fixed at $0$.) This operator 
is treated in section 5 using several considerations including Lemma 4.1, (4.9), resolvent
equation and semi-algebraic sets. The decay rate $\alpha_{i+1}$ depends on $\gamma_{M_0}$,
where $M_0\ll N=M^{i+1}$ is determined in (5.5), cf. also (6.1-6.8).
So (4.16) prevents the deterioration of $\alpha_i$ as $i\to\infty$ and we will have $\alpha_i>\alpha/2$
for all $i$ in section 6.

We prove Lemma 4.1 using iteration.
The two pillars of this iteration are semi-algebraic set techniques and Cartan-type theorem for analytic matrix 
valued functions (see [B5, BGS]).
The general construction of the iteration is the same as in our previous paper [BW]. 

\noindent 
{\it The initial estimate $(0^{\text{th}}$ step)}

In view of (4.4) and the Newton scheme (3.9), choose $M>2p$, such that
$$
\Cal S\cup -\Cal S \subset [-M, M]^{d+\nu}.\tag 4.17
$$
\proclaim{Lemma 4.2} Fix $0<\sigma<\beta<1$, there exists $M_0=M_0(d+\nu,\sigma,\beta)$ 
such that for all $M\geq M_0$,  
there exist $\epsilon_0$, $\delta_0$ such that 
$$\aligned 
&\Vert [T_M^\theta]^{-1} \Vert < e^{M^\beta}\qquad (0<\beta<1),\\
&|[T_M^\theta]^{-1}(\ell, \ell')|<e^{-\alpha|\ell-\ell'|}\qquad (\alpha \text{ as in }(H3)),\\
&|\ell-\ell'|>M/10, M\geq M_0,\endaligned$$  for $0<\epsilon\leq\epsilon_0$, $0<\delta\leq\delta_0$, 
$\theta\in\Bbb R\backslash \Cal B_{\beta,\alpha}(M)$, $\text{mes }\Cal B_{\beta,\alpha}(M)\leq e^{-M^\sigma}$, all $x\in\Bbb R^{\Bbb Z^d}\backslash\Bbb R^{\nu}$, 
all $\omega\in [0,1)^\nu$. (Recall from (4.3-4.5) the $x$ and $\omega$ dependence of $T_N^\theta(y'_i)$.)
\endproclaim
\demo {Proof}
We use Neumann series in $\epsilon, \delta$ to estimate $T^{-1}_{M}$.
We require
$$
\cases
|\theta +n\cdot\omega+v_j|> 2\max (e^{-M^\beta}, (\epsilon+\delta)^{1/2}),\\
|-\theta-n\cdot \omega+v_j|> 2\max (e^{-M^\beta}, (\epsilon+\delta)^{1/2}),\endcases
\forall (j, n) \in [-M, M]^{d+\nu} \backslash \{\Cal S\cup- \Cal S \}.\tag 4.18
$$
Clearly, (4.18) holds away from a set of $\theta \in\Bbb R$ of measure at most
$$
4(2M+1)^{d+\nu}\max (e^{-M^\beta}, (\epsilon+\delta)^{1/2}).\tag 4.19
$$

Choose $\epsilon$, $\delta$ such that $(\epsilon+\delta)^{1/2}\leq e^{-M^\beta}$, which can
be satisfied if $0<\epsilon\leq\epsilon_0$, $0<\delta\leq\delta_0$ with 
$\epsilon_0=\delta_0=\frac{1}{2}e^{-2M^\beta}$. From (4.19), we need 
$4(2M+1)^{d+\nu}e^{-M^\beta}\leq e^{-M^\sigma}$ ($0<\sigma<\beta<1$). This leads to
$M\geq M_0(d+\nu,\sigma,\beta)$. 

On the complement of the set defined in (4.18), using Neumann series in 
$\epsilon, \delta$ for $T_M^{-1}$ and (H3), we verify
that
$$\aligned 
&\Vert [T_M^\theta]^{-1} \Vert < e^{M^\beta}\qquad (0<\beta<1),\\
&|[T_M^\theta]^{-1}(\ell, \ell')|<e^{-\alpha|\ell-\ell'|}\qquad (\alpha \text{ as in }(H3)),\\
&|\ell-\ell'|>M/10, M\geq M_0,\endaligned$$ 
for $0<\epsilon\leq\epsilon_0$, $0<\delta\leq\delta_0$.
The probability set at this scale $X_M$ and the frequency set at this scale $\Omega_M$,
on which and on the complement of which (4.14) holds satisfy
$$\text{\rm mes\,} X_M=1, \text{\rm mes\,}\Omega_M=0\tag 4.20$$
\hfill $\square$
\enddemo

This direct perturbation argument is the same as in \cite{BW}. Note that (4.20) entails that (4.14) holds for all
$\omega, v_j$, as in the $0^{\text{th}}$ step, the invertibility is  entirely provided by shifting in $\theta$.
There is {\it no} bad site. We will only use the above lemma for the initial set of scales.

\noindent
{\it The iteration}

We now prove Lemma 4.1 using iteration from scale $N_0$ to $N_0^C=N$. ($C$ assumed large).
We call an elementary region $\Lambda(N_0)$ {\it bad}, if the first $2$ inequalities in (4.9) are violated.
As in [BW], we need to perform an incision in the frequency space, in order that inside any $\Lambda(N)\in\Cal E\Cal
R(N)$, there are at most $N^{1-}$ pair-wise disjoint bad elementary regions at scale $N_0$, 
where  $N^{1-}$ means $N^a$ ($0<a<1$). For
technical reasons (cf. \cite{BGS, BW}), this requirement pertains to all elementary regions $\Lambda(N')$, $N_0\leq N'\leq
2N_0$ and not simply at $N'=N_0$. 
For later constructions in section 5, it is important to note once again that this set is independent of
$\Cal V=\{v_{j_k}\}^\nu_{k=1}$.

Let $\Lambda(N)\in \Cal E\Cal R(N), \bar \Lambda(N)$ be its projection onto $\Bbb Z^d$;
$$
\Cal T{(N)}= \{[-N_0, N_0]^d\times [-N, N]^\nu\}\cap \Lambda(N),\tag 4.21
$$
$\bar{\Cal T}(N)$ be its projection onto $\Bbb Z^d$.
Denote by $\overline{\Cal E\Cal R}(N)$ the projection of $\Cal E\Cal R(N)$ onto $\Bbb Z^d$.
Note that $\bar\Lambda(N)\in\overline{\Cal E\Cal R}(N)$ can be of much smaller diameter than $2N$.
\smallskip
\noindent {\it Number of bad elementary regions at scale $N_0$ disjoint from ${\Cal T}(N)$.}

Using (H3) in (4.2-4.4), the region $\Lambda(N)\backslash \Cal T(N)$ can be treated perturbatively.
We make separate incisions in the probability space  and the frequency space. We first make 
incisions in the probability space. Toward that end, we look at 
$$\
\Lambda(N')\in\Cal E\Cal R(N'),\, \Lambda(N')\subset\Lambda(N), \, \Lambda(N')\cap {\Cal T}(N)=\emptyset\quad
(N_0\leq N'\leq 2N_0).\tag 4.22
$$
Let 
$$
X_N'\subset \Bbb R^{\{\Bbb Z\backslash [-N_0, N_0]\}^d}\tag 4.23
$$
be the probability set such that there is at most {\bf 1} (pair-wise disjoint) $\bar\Lambda(N')$ ($N_0\leq
N'\leq 2N_0$), satisfying :
$$
\cases
\bar\Lambda(N')\in \overline{\Cal E\Cal R}(N'),\\
\bar\Lambda(N')\subset[-N, N]^d,\\
\bar\Lambda(N')\cap [-N_0, N_0]^d=\emptyset,
\endcases
\tag 4.24
$$
where $(A1) $ is violated for some $E\in I$, $I=\sigma(H)$, the set defined in (1.6).

Using Theorem A, with $m=\gamma$,
$$
\align
\text{mes\,} X_N' &> \bigg(1-\frac 1{N_0^{2p'}}\bigg)^{\Cal O(1)(N^d\times N_0^{d+1})^2}\\
&\geq 1-\frac{\Cal O(1)}{N^{\frac 2C(p'-d(C+1)-1)}}\quad \big(p'>d(C+1)+1\big),\tag 4.25
\endalign
$$
where $\Cal O(1)$ is a universal geometric constant; $p'$ is to be determined at the conclusion 
of the proof of Lemma 4.1. We used $N=N_0^C$,
the second factor
$N_0^{d+1}$ comes from the estimate on the number of elementary regions of sizes 1 to $2N_0$ associated to each
lattice site:
$\Cal O{(1)}\sum_{N_0'=1}^{2N_0} N_0'{}^d$ \big(cf. (4.11)\big) the exponent is an upperbound on the number of
pairs of elementary regions of sizes up to $2N_0$ in $\Lambda(N)$.

Given two elementary regions $\Lambda_1(N')$, $\Lambda_2(N')$, we say that they are convex-disjoint if their
convex envelops are disjoint.
(This is in order that we have (4.23-4.25) at our disposal.) To control the number of bad 
elementary regions at scale $N_0$, we now make additional incisions in the frequency space. Recall
that (4.23, 4.24) pertain only to the projected elementary regions in $\Bbb Z^d$.

We are now ready to prove

\proclaim{Lemma 4.3} Fix $x\in X_N'\cap \tilde X_{N_0}$, where $X'_N$ is the set defined in (4.23, 4.24)
and $\tilde X_{N_0}$ is defined as in (4.23, 4.24), but with $N_0^{1/C}$ replacing $N_0$. 
There exists a set $\Omega_N'$,
$$
\text{\rm mes\,} \Omega_N' \leq e^{-N^{\frac{\beta}{2C}}},\tag 4.26
$$
such that if $\omega\not\in \Omega_N'$, then for any $\Lambda(N)\in\Cal E\Cal R(N)$
{\text any} fixed $\theta$, there are at most {\bf 2} convex-
disjoint
bad $\Lambda(N')\in\Cal E\Cal R(N')$, $\Lambda(N')\cap \Cal T(N)=\emptyset$, $N_0\leq N'\leq 2N_0$ in $\Lambda(N)$
($N=N_0^C$). Moreover $\Omega_N'$ is semi-algebraic with degree bounded above by $\Cal O(1)N^{6d+\nu}$ and is contained in the union of at most $\Cal O(1) N^{6d+\nu}$ components.
\endproclaim
\noindent
{\it Remark.} $\Omega_N'$ is independent of $\Cal V$, $y_{k}$. Observe also that we need localization
information on the random Schr\"odinger operators at $2$ scales, $N_0$ and $N_0^{1/C}$.

\noindent
{\bf Proof.} In view of (4.3, 4.4, 4.7, H3), for $\Lambda(N_0)$ such that $\Lambda(N_0)\cap \Cal T(N)=\emptyset$,
$\delta S$ can be treated as a small perturbation.
We only need to ensure the invertibility of $D^\theta_{\Lambda(N_0)}$.
Assume $\Lambda(N_0), \Lambda'(N_0)$, and $\Lambda''(N_0)$ are 3 convex-disjoint {\it bad} elementary regions.
So $\exists (n, j)\in \Lambda(N_0), (n', j')\in \Lambda'(N_0)$, $(n'', j'')\in \Lambda''(N_0)$, such that
$$
\align
|\theta+n\cdot\omega+\mu_j|<2e^{-N_0^\beta}\, \text { or }\, |-\theta-n\cdot \omega+\mu_j|< 2e^{-N_0^\beta};\\
|\theta+n'\cdot \omega+\mu_{j'}|<2e^{-N_0^\beta}\, \text { or }\, |-\theta-n'\cdot \omega+\mu_{j'}|<2
e^{-N_0^\beta};\tag 4.27
\endalign
$$
and 
$$
|\theta+n''\cdot \omega+\mu_{j''}|<2e^{-N_0^\beta}\text { or } |-\theta-n''\cdot\omega+\mu_{j''}|<2e^{-N_0^\beta},
$$
where $\mu_j$, $\mu_{j'}, \mu_{j''}$ are eigenvalues of $\bar\Lambda(N_0), \bar{\Lambda'}(N_0)$, $\bar \Lambda'' (N_0)$
respectively.

(4.27) implies that $\exists m, \lambda$ such that
$$
|m\cdot \omega+\lambda|< 4e^{-N_0^\beta},\tag 4.28
$$
where $m=\pm(n-n'), \text { or } \pm(n'-n'') \text { or } \pm (n-n'')$,
$$
\lambda=\mu_j-\mu_{j'}, \text { or } \mu_{j'}-\mu_{j''}\text{ or } \mu_{j''}-\mu_{j}.\tag 4.29
$$
We use the same argument as in the proof of Lemma 2.3 of [BW], which we summarize below.

There are two possibilities: $m=0, m\not= 0$.
When $m=0$, from pairwise disjointness (4.24, 4.25), (A6) implies
$$
|\lambda|>e^{-N_0^{\beta'}} \qquad (0<\beta'< \beta)\tag 4.30
$$
which contradicts (4.28).

When $m\not= 0$, (4.28) corresponds to at most
$$
\Cal O(1) N^\nu\cdot (N^d\cdot N_0^{2d+1})^2<\Cal O(1)N^{6d+\nu}\tag 4.31
$$
number of equations. Since each equation in (4.28) is a monomial of degree $1$ in $\omega$,
the excluded set $\Omega'_N$ is of degree less than $\Cal O(1)N^{6d+\nu}$. 
Since $|\omega|$ may be assumed to be bounded for each such equation, we exclude a set 
in $\omega$  of measure $\Cal O(1) e^{-N_0^\beta}$.
It is simple to see that for each equation
in (4.28), the excluded set in $\omega$ has $1$ single component. 

So in conclusion, for fixed $x\in X'_N$, $\omega\in{\Omega_N'}^c$, there are at most $2$ 
pair-wise disjoint $\Lambda(N')\in\Cal {ER}(N')$,  $\Lambda(N')\cap\Cal T(N)=\emptyset$
such that the first inequality in (4.9) is violated by using (H3) and a simple perturbative 
argument. 

Assume $\Lambda(N')$ is such that the first inequality of (4.9) is satisfied. So
$|\pm\theta\pm n\cdot\omega+\mu_j|\geq 2e^{-N_0^\beta}$ for $n,\,j\in\Lambda(N')$ 
from above considerations. To obtain the second inequality we proceed as follows. Let
$\bar\Lambda(N')$ be the projection of $\Lambda(N')$ onto $\Bbb Z^d$. In view of the 
restriction in the third expression in (4.9), we may assume $\text{diam }\bar\Lambda(N')\geq\frac{N_0}{10}$. We cover $\bar\Lambda(N')$ with elementary regions $\bar\Lambda(N_0^{1/C})$
of diameter $2N_0^{1/C}$.

Since $x\in X'_N\cap\tilde X_{N_0}$ and on $\tilde X_{N_0}$ for all $E$ (of the form 
$E=\theta+n\cdot\omega$ here), there is at most $1$ (pair-wise disjoint) $\bar\Lambda(N_0^{1/C})$
in $\bar\Lambda(N')$ where (A1) (with $m=\gamma$) is violated. Using the resolvent equation,
(A1) and the estimates  $|\pm\theta\pm n\cdot\omega+\mu_j|\geq 2e^{-N_0^\beta}$ for $n,\,j\in\Lambda(N')$ for the bad $\bar\Lambda(N_0^{1/C})$, we obtain exponential decay in the $j$ direction
for $\bar\Lambda(N')$ for all $E$ (of the form $E=\theta+n\cdot\omega$ here). 

We obtain the second inequality in (4.9) for this $\Lambda(N')$ by another application of 
resolvent series in $\delta S$ and using (H3). This holds for all $\Lambda(N')$ such that
the first inequality of (4.9) is satisfied, in view of the definition of $\tilde X_{N_0}$. 
Using $N_0=N^{1/C}$, we obtain the lemma.
$\hfill\square$

\noindent{\it Number of bad elementary regions at scale $N_0$ intersecting ${\Cal T}(N)$.}

We now estimate the number of bad $\Lambda(N')\in\Cal E\Cal R(N')$,
$\Lambda(N')\subset \Lambda(N)$, $\Lambda(N')\cap \Cal T(N)\not= \emptyset$, $\Cal T(N)$ as in (4.21), $N_0\leq
N'\leq 2N_0$, using semi-algebraic set techniques. Here it is important to emphasize the $\Bbb Z^d$
coordinate of the center of elementary regions as the linearized operator is not a Toeplitz operator
in the $\Bbb Z^d$ variable. We look at elementary regions with centers in $\{0\}\times\Bbb Z^d$. We write
$\Cal E\Cal R(N', j)$ for the set of elementary regions centered at $j\in\Bbb Z^d$.
For any $\Lambda(N', j)\in \Cal E\Cal R(N', j)$,
let $\Cal B(N', j)\overset\text{def }\to=\Cal B_{\beta,\gamma}(N', j)$ be a set such that on $\Cal B^c(N', j)$,
(4.9) hold. (Later for more general elementary 
regions centered in $i\in\Bbb Z^{d+\nu}$, we will use the same notations as used here.)

Assume $\exists X_{N', j}, \Omega_{N', j}$ such that for $x\in X_{N', j}$, $\omega\in DC_{A, c}(2N)\backslash
\Omega_{N', j}$,
$$
\text{mes\,} \Cal B(N', j)\leq e^{-N_0^\sigma}\quad(N_0\leq N'\leq 2N_0,\quad 0<\sigma<1).\tag 4.32
$$
Let 
$$
\align
&X_{N_0}'' {\overset\text{def }\to=} \, \bigcap_{j\in[-3N_0, 3N_0]^d} \,\bigcap_{N_0\leq N'\leq 2N_0}\, X_{N',
j},\\ &\Omega_{N_0}'' {\overset\text{def }\to =} \, \bigcup_{j\in[-3N_0, 3N_0]^d} \,\bigcup_{N_0\leq N'\leq 2N_0}\,
\Omega_{N', j},\\ 
&\Cal A {\overset\text{def }\to =}\, \bigcup_{j\in[-3N_0, 3N_0]^d}\, \bigcup_{N_0\leq N'\leq
2N_0}\,\bigcup_{\Cal E\Cal R(N',j)}\, \Cal B(N', j),\\
&\text{mes }\Cal A\leq \Cal O(1)N_0^{2d+\nu} e^{-N_0^\sigma}\quad \text {from } (4.32).\tag 4.33\endalign$$

\proclaim
{Lemma 4.4}
Let $N=N^C_0$.
For any fixed $\theta\in \Bbb R$, let
$$
I=\{n\in [-N, N]^\nu \, |n\cdot\omega+\theta\in \Cal A\}.\tag 4.34
$$
Fix $x\in X_{N_0}''$.
Then for $\omega\in DC_{A, c} (2N)\backslash \Omega_{N_0}''$, 
$$
|I|\leq \Cal O_{d, v}(1) N_0^{6(d+\nu)} =\Cal O_{d, v}(1)N^{\frac{6(d+\nu)}{C}}{\overset\text{def }\to
=}N^{1-b_0},\tag 4.35
$$
$0<b_0<1$ and assuming $6(d+\nu)<C<N_0^{\sigma/2}$. Hence there are at most $\Cal O(1)N^{1-b_0}$
($\Cal O(1)$ a universal geometric constant) pair-wise disjoint bad $\Lambda(N')\in\Cal {ER}(N')$,
$\Lambda(N')\cap\Cal T(N)\neq\emptyset$, $N_0\leq N'\leq 2N_0$ in $\Lambda(N)$ ($N=N_0^C$).
\endproclaim

\noindent
{\bf Proof.}
Since the Green's function is the ratio of two determinants and the norm of the Green's function can be replaced by
its Hilbert-Schmidt norm, (4.9) can be reexpressed as polynomial inequalities in $\theta$.
$\Cal A$ is therefore a {\it semi-algebraic} set.
(See [Ba, section 7 of BGS].)

$\Cal A$ is defined by
$$
\align
&\Cal O(1)N_0^d\cdot N_0^{d+\nu}\cdot N_0^{2(d+\nu)}\\
&\cong\Cal O(1)N_0^{4d+3\nu}\tag 4.36
\endalign
$$
number of polynomials,
$N_0^d$ for number of centers, $N_0^{d+\nu}$ number of elementary regions per center, $N_0^{2(d+\nu)}$ number of 
matrix elements.
Each polynomial is of degree $N_0^{2(d+\nu)}$ in $\theta$ (as one squares the matrix elements.)

Basu's Theorem [Ba], restated as Theorem 7.3 in [BGS], then gives that the number of connected components in $\Cal
A$ does not exceed $N_0^{6(d+\nu)}$.
If there are $n, n', n\not= n'$, such that $n, n'$ belong to the same connected component of $\Cal A$, then from
the last inequality in (4.33),
$$
|(n-n')\cdot\omega|\leq \Cal O(1) N_0^{2d+\nu} e^{-N_0^\sigma}.\tag 4.37
$$
Since $n, n'\in [-N, N]^\nu, n-n'\in [-2N, 2N]^\nu, N=N_0^C$, $\omega \in DC_{A, c}(2N)$ is in contradiction with
(4.37) for $C<N_0^{\sigma/2}<\frac{N_0^{\sigma}}{2A\log N_0}$ (assuming $N_0\gg 1$), so there can be at most 1
integral point in a connected component of $\Cal A$. We therefore obtain (4.35). 

Let $i=(j,n)\in\Bbb Z^{d+\nu}$. Since $\Lambda(N',i)=\Lambda(N',j)+n$ and $T$ is a Toeplitz operator
in the $\Bbb Z^\nu$ variable, we obtain the second conclusion of the lemma.
$\hfill\square$

\noindent{\it Remark.} $C$ will be a {\it fixed} expansion factor. So the upper bound on $C$ will 
be satisfied for all $N\geq N_0$ as soon as it is satisfied for an initial $N_0$.
\smallskip
\noindent
{\it A large deviation estimate in $\theta$.}

Lemmas 4.3, 4.4 combined give that the number of bad elementary regions at scale $N_0$ in $\Lambda (N)$ is
at most $N^{1-b_0}$, where 
$$N=N_0^C,\quad b_0=1-\frac{6(d+\nu)}{C}\quad (6(d+\nu)<C<N_0^{\sigma/2})$$
from (4.35). This  
enables us to use a Cartan-type theorem for analytic matrix valued functions (see [B5]) to prove a
large deviation estimate on $\Vert T_N^\theta(y_{i'})^{-1}\Vert$, necessary for the proof of Lemma 4.1.
The proof of the lemma is very similar to the one in [BW] (See also [BGS]), after using (4.15) to appropriately
adjust
$y_{i'}$ according to the scale $N$.
So we state (without details of the proof)

\proclaim
{Lemma 4.5}
Let $b_0$, $\beta$, $\sigma$, $\gamma$ be fixed positive numbers, so that
$$
0<b_0, \beta, \sigma<1 \text { and } \beta+b_0> 1+3\sigma.\tag 4.38
$$
Let $N_0\leq N_1$ be positive integers satisfying
$$
\bar N_0(\beta, \sigma, \gamma)\leq 100 N_0\leq N_1^\sigma\tag 4.39
$$
with some large constant $\bar N_0$ depending only on $\beta, \sigma, \gamma$.
Assume that for any $N_0\leq L\leq N_1$, any $\Lambda(L)\in\Cal E\Cal R(L, i)$, $i\in\Bbb Z^{d+\nu}$,
$$
\text{\rm mes\,} \Cal B_{\beta,\gamma}(L, i)\leq e^{-L^\sigma}.\tag 4.40
$$
Let $\bar X_N$, $\bar\Omega_N$ be the sets such that on $\bar X_N$ and $DC_{A,c}(2N)\backslash\bar\Omega_N$,
(4.40) holds for all $i\in [-N,N]^{d+\nu}$, $L\in [N_0,N_1]$.
Let $X_N'$, $\tilde X_{N^{1/C}}$ be the sets defined in (4.23, 4.24) and Lemma 4.3, $\Omega_N'$ the corresponding frequency set as in Lemma 4.3.
If
$$
\aligned
&x\in X_N' \cap \tilde X_{N^{1/C}}\cap \bar X_N,\\
&\omega\in DC_{A, c} (2N)\backslash\{ \Omega_N'\cup \bar\Omega_N\},\\
\endaligned
\tag 4.41
$$
then
$$
\text {\rm mes\,} \big\{\theta| \, \Vert [T_N^\theta]^{-1}\Vert> e^{N^\beta}\big\} < e^{-N^{3\sigma}},\tag 4.42
$$
where $T_N$ is restriction of $T$ to any $\Lambda\in\Cal E\Cal R(N)$, the elementary regions centered at $0$, provided
$N_0^{C_1} \leq N\leq N_1^{\sigma C_1}$, $C_1\gg\max (\frac 1\sigma, 6(d+\nu))$
 depending only on $\beta,\sigma$.
\endproclaim

\noindent{\it Remark.} We note that from (4.23, 4.24), the dependence of the probability set $X_{L, i}$ (on 
which (4.40) holds) on $i\in\Bbb Z^{d+\nu}$ is only through the $\Bbb Z^d$ coordinate. For simplicity, we keep
the notation $X_{L, i}$. The set $\Omega_{L, i}$ (on the complement of which (4.40) holds), on the other hand, does have full dependence on $i$.

\noindent {\it A summary of the proof.}

We use analytic and harmonic function theory together with a 2-scale (in the range $[N_0, N_1]$) analysis to 
control the measure of the set in (4.42) at scale $N\gg N_0$. (See the proofs of Lemma 4.4 in \cite{BGS}
and Lemma 6.2 in \cite {BW}.) For a given $N$, let these $2$ scales $L_1$, $L_2\in [N_0, N_1]$, $L_1<L_2$,
satisfy $$\log L_1\sim\frac{1}{C_1}\log N,\quad \log L_2\sim\frac{1}{\sigma}\log L_1\sim \frac{1}{C_1\sigma} \log
N,$$
$C_1$ as in the lemma, $C_1\gg 1/\sigma$. 
We reiterate
the main line of arguments below.

\itemitem {$\bullet$} Fix $\theta$. Let
$$
\align
\Lambda_{*}^c= \{m\in\Lambda(N)|\exists\, \Lambda_1 \in \Cal ER(L_1), N_0&\leq L_1\leq 2N_0,
\Lambda_1\subset m+[-L_1,L_1]^{d+\nu}, \Lambda_1 \text { is bad}\},\tag 4.43
\endalign
$$
For $x$, $\omega$ satisfying (4.41), 
$$|\Lambda_{*}^c|\leq N^{1-b_0} \quad (b_0>0),\tag 4.44$$
by Lemmas 4.3, 4.4. Since $X''_{N_0}\supset\bar X_N$, $\Omega''_{N_0}\subset\bar\Omega_N$, where
$X''_{N_0}$, $\Omega''_{N_0}$ as defined by the first two equations of (4.33).

\itemitem {$\bullet$} Let $\Lambda_{*}$ be, roughly 
speaking, the complement of the set in (4.43). For more precise definition, which requires a partition
of $\Lambda$, see the beginning of the proof of Lemma 4.4 in \cite{BGS}. Using an elementary resolvent
expansion, Lemma B in the appendix, which is Lemma 2.2 of \cite{BGS} reiterated, we obtain an upper bound on 
$\Vert[T_{\Lambda_{*}}^{\theta}]^{-1}\Vert$ at fixed $\theta$
by using the decay estimate on the $\Lambda_1$'s, elementary
regions at scale $L_1$, in $\Lambda_{*}$. By  definition they are all {\it good}. By standard Neumann series arguments,
this bound is preserved inside the disk $B(\theta,e^{-N_0})\subset\Bbb C$. 

\noindent{\it Remark.} We have control over the size of $\Lambda_{*}$ via (4.43, 4.44), but not its geometry.
Typically
$\Lambda_{*}$ is {\it non-convex}, hence the need for elementary regions which are more general than cubes, in
particular L-shaped regions, in view of Lemma B. 

\itemitem {$\bullet$} Define a matrix-valued analytic function $A(\theta')$ on $B(\theta, e^{-N_0})$ as
$$\aligned A(\theta')&=R_{\Lambda_{*}^c} T_N^{\theta'} R_{\Lambda_{*}^c}\\
&\quad -R_{\Lambda_{*}^c} T_N^{\theta'} R_{\Lambda_{*}}{[T_{\Lambda_{*}}^{\theta'}]}^{-1}
R_{\Lambda_{*}} T_N^{\theta'} R_{\Lambda_{*}^c}\endaligned\tag 4.45$$
where $\Lambda_{*}^c=\Lambda\backslash \Lambda_{*}$, $R_{\Lambda_{*}}$, $R_{\Lambda_{*}^c}$ are projections.
From (4.44), $A(\theta')$ is a rank $\Cal O(N^{1-b_0})\times \Cal O(N^{1-b_0})$ matrix. The
raison d'etre of introducing $A(\theta')$ is the following inequality:
$$\Vert A(\theta')^{-1}\Vert\lesssim \Vert (T_N^{\theta'})^{-1}\Vert \lesssim e^{2N_0}\Vert
A(\theta')^{-1}\Vert,$$ (see Lemma 4.8 of \cite{BGS}). So to bound $\Vert (T_N^{\theta'})^{-1}\Vert$, 
it is sufficient to bound
$\Vert A(\theta')^{-1}\Vert$, which is of smaller dimension.

\itemitem {$\bullet$} Toward that end, we introduce an intermediate scale $L_2$, $\log L_2\sim (\log
L_1)/\sigma>\log L_1$. We work in an interval $\Theta=\{\theta'||\theta'-\theta|<e^{-N_0}\}$. Using (4.40) for the
$\Lambda_2$'s at scale $L_2$ and in
$\Lambda$, Lemma B, we obtain an upper bound on $\Vert [T_N^{\theta'}]^{-1}\Vert$
except for a set of $\theta'$ of measure smaller than $e^{-\Cal O(L_2^\sigma)}$. So there exists $y\in\Theta$, such
that we have {\it both} an upper bound on $\Vert A^{-1}(\theta')\Vert$ at $\theta'=y$, hence a lower bound on
the smallest eigenvalue of $A(\theta)$ and an
a priori upper bound on $\Vert A(\theta')\Vert$, which comes from the boundedness of $T_N^{\theta'}$
and the bound on $[T_{\Lambda^*}^{\theta'}]^{-1}$ (see (4.45)).

\itemitem {$\bullet$} Transfering the estimates on $\Vert A(\theta')\Vert$, $\Vert A^{-1}(\theta')\Vert$ into
estimates on $\log|\det A(\theta')|$, which is subharmonic and using either Cartan-type theorem (see sect. 11.2
in \cite {Le}) or proceeding as in the proof of Lemma 4.4 of \cite{BGS} or Chap XIV of \cite{B5}, we obtain the
lemma by covering the interval
$I=(-\Cal O(N_1^{\sigma C_1}),\Cal O(N_1^{\sigma C_1}))$ with intervals of size $e^{-N_0}$. (Recall (4.39) and
that for all $\theta\notin I$, $T_N^{\theta'}$ is automatically invertible.)
$\hfill\square$

\noindent
{\it Iteration Lemma}

To obtain exponential decay of $T_N^{-1}$, we need

\proclaim
{Lemma 4.6} Suppose $M, N\in\Bbb N^+$ are such that for some $0<\tau <1$
$$
N^\tau\leq M\leq 2N^\tau.\tag 4.46
$$
Let $\Lambda_0\in \Cal E\Cal R(N)$ be an elementary region with the property that for all $\Lambda \subset\Lambda_0,
\Lambda \in\Cal E\Cal R(L)$ with $M\leq L\leq N$,
$$
\Vert (T_L^\theta)^{-1} \Vert\leq e^{L^\beta} \qquad (0<\beta<1).\tag 4.47
$$
We say that $\Lambda \in\Cal E\Cal R(L), \Lambda\subset\Lambda_0$ is good, if in addition to (4.47), 
$$
\Vert (T_L^{\theta})^{-1} (k, k')\Vert\leq e^{-\gamma|k-k'|}\tag 4.48
$$
for all $k, k'\in\Lambda, |k-k'|>L/10$.
Otherwise $\Lambda$ is called bad.

Assume that for any family $\Cal F$ of pairwise disjoint bad $M'$-regions in $\Lambda_0$ with
$$
\align
M+1\leq M'&\leq 2M+1\\
\sharp\Cal F&\leq N^\beta (0<\beta<1).\tag 4.49
\endalign
$$
Then one has
$$
|({T_N^{\theta}})^{-1} (k, k')|\leq e^{-\gamma'|k-k'|}\tag 4.50
$$
for all $k, k'\in\Lambda_0, |k-k'|> N/10$, $\gamma'=\gamma-N^{-\delta'} (\delta'>0)$, provided $N\gg N_0(\beta,
\tau, \gamma)$.
\endproclaim

The proof of the above lemma is written out in detail in [BGS]. The only difference is that instead of
being tridiagonal, $T_L^\theta$ has exponentially decaying off-diagonal elements from (H3). So 
$\gamma=\gamma(\alpha)$. The proof goes through.
So we do not repeat it here. (See also [BW].) The gist is as follows.

The exponential decay estimate at scale $N$ in (4.50) is obtained from exponential decay estimate in (4.48)
at smaller scales $M'$ by using (4.15), the norm estimate in (4.47) and the resolvent identity. To implement this,
we use a sequence of scales $M_{j+1}=M_j^{C'}$, with $M_0=M$, $C'$ is chosen such that $C'\beta<1$, $C'\tau\leq 1$.
For each elementary region at scale 
$M_{j+1}$, $\Lambda(M_{j+1})$ and for each $k\in \Lambda(M_{j+1})$,
we exhaust
$\Lambda(M_{j+1})$ by an increasing sequence of annuli centered at $k$ of width $2M_j$, or more precisely the
intersection of this sequence with $\Lambda(M_{j+1})$. Roughly speaking, an annulus is good, if it does not
intersect a bad cube of the previous scale $M_j$.

In each of the connected components of the complement of the bad annuli, we apply the resolvent identity
using the estimate in (4.48) for elementary regions of size $M_{j}$. In the bad annuli, we use 
(4.47). From (4.49), the number of bad annuli is at most sublinear in $M_{j+1}$. Using a multiscale
induction argument to reach the scale $N$, we obtain the exponential decay in (4.50), when $|k-k'|>N/10$. 
$\delta'$ is determined from (4.47, 4.49), $\delta'\simeq\tau(1-\beta C')$.
$\hfill\square$

Let $\Cal B_{\beta, \gamma}(N, i)$, $i\in\Bbb Z^{d+\nu}$ be a set such that on the complement, (4.9) hold.
When $i=0$, we write $\Cal B_{\beta, \gamma}(N, 0)=\Cal B_{\beta, \gamma}(N)$.
As before let 
$$\aligned&X_{N,i},\,\Omega_{N,i}\text{ be the probability, frequency subsets on which and}\\
&\text {on the complement of which
mes }\Cal B_{\beta,\gamma}(N,i)\leq e^{-N^\sigma}\endaligned\tag 4.51$$
Combining Lemmas 4.5, 4.6 with (4.23-4.25), Lemmas 4.3, 4.4, we obtain

\proclaim
{Lemma 4.7} Assume that for any  $\tilde N_0\leq N_0\leq \tilde N_0^C, 
\,\max (\frac 1\sigma, 6(d+\nu)\ll C<\tilde
N_0^{\sigma/2}$ ($\tilde
N_0(\epsilon, \delta)\gg 1$), any $\Lambda(N_0)\in\Cal{ER}(N_0, i)$, $i\in\Bbb Z^{d+\nu}$ 
$$
\text{\rm mes\,} \Cal B_{\beta, \gamma}(N_0, i) \leq e^{-N_0^\sigma} \qquad (0<\sigma<1).\tag 4.52
$$
Let $[\tilde N_0^C, \tilde N_0^{C^2\sigma}]$ be the next interval of scales. For any $N\in [\tilde N_0^C, \tilde
N_0^{C^2\sigma}]$, write $N=N^C_0$ with $N_0\in [\tilde N_0, \tilde N_0^{C}]$.
Let $X'_{N}$, $\tilde X_{N^{1/C}}$ be the sets defined in (4.23, 4.24) and Lemma 4.3, $\Omega'_{N}$ the set defined
in Lemma 4.3, satisfying 
$$
\aligned
&\text { \rm mes\,} X_N'\geq 1-\frac{\Cal O(1)}{N^{{\frac 2C}(p'-d(C+1)-1)}}\qquad (p'>d(C^2+1)+1),\\
&\text {\rm mes\,} \tilde X_{N^{1/C}}\geq 1-\frac{\Cal O(1)}{N^{{\frac 2{C^2}}(p'-d(C^2+1)-1)}}\qquad (p'>d(C^2+1)+1),\\
&\text{\rm mes\,} \Omega_N'\leq e^{-N^{\frac {\beta}{2C}}},
\endaligned\tag 4.53
$$
in view of (4.25, 4.26). 

Let $\bar X_N$, $\bar\Omega_N$ be the sets such that on $\bar X_N$ and $DC_{A, c}(2N)\backslash\bar\Omega_N$,
(4.52) holds for all $i\in [-\tilde N_0^{C^2\sigma}, \tilde N_0^{C^2\sigma}]^{d+\nu}$, all $N_0\in [\tilde N_0,
\tilde N_0^C]$. If
$$
\align
&x\in X_N' \cap \tilde X_{N^{1/C}}\cap \bar X_N{\overset\text{\rm def }\to=} X_N,\\
&\text{mes }X_N\geq 1-1/N^{p''} \text{ with }p''=p''(p',C,d+\nu)>1 \text{ by choosing }p'\text{ large enough};\\
&\omega\in DC_{A, c} (2N)\backslash\{ \Omega_N'\cup\bar\Omega_N\}
{\overset\text{def }\to=} DC_{A, c} (2N)\backslash \Omega_N,\tag 4.54\\
&\text{mes }\Omega_N\leq e^{-N^{\frac{\beta}{3C^2}}},\\
&\Omega_N \text{ is semi-algebraic with degree less than } N^{C^2(d+\nu)}
\endalign
$$
then
$$
\text{\rm mes\,} \Cal B_{\beta, \gamma'}(N)\leq e^{-N^\sigma}\qquad (0<\sigma< 1).\tag 4.55
$$
where $\gamma'=\gamma-N^{-\delta'}$ ($\delta'>0$).
\endproclaim

\demo{Proof}
Applying Lemma 4.5 using (4.52), we obtain the large deviation estimate on $\Vert
[T_N^\theta]^{-1}\Vert$. Choosing $0<\tau<\frac{1}{C\sigma}$, for $x$, $\omega$ in the sets defined in
(4.54), (4.49) is satisfied. (Here we need the definition that an $N^{\tau}$-region is bad if it intersects
a bad $N_0$-region. Otherwise it is good. On the good $N^{\tau}$-region, (4.48) is obtained by using resolvent
expansion (Lemma B) and (4.9) for  $N_0$-regions.) Hence Lemma 4.6 is available and we obtain (4.55).
The estimates on $X_N$, $\Omega_N$ follow from (4.51, 4.53) and the constructions in Lemme 4.3-4.5.
$\hfill\square$
\enddemo 
\smallskip
\noindent
{\bf Proof of Lemma 4.1.}
Assume $\bar N_0$ (to be determined below) is such that Lemma 4.7 is available. For the scales
$\bar N_0\leq N\leq\bar N_0^C$, we use Neumann series in $\epsilon$ and $\delta$ a la Lemma 4.2 and its proof.
For the scales $N\geq \bar N_0^C$, we use Lemma 4.7.

From Lemma 4.2, we need $$\bar N_0\geq M_0(d+\nu,\sigma,\beta)\tag 4.56$$
and $$\epsilon_0=\delta_0=\frac{1}{2}e^{-2\bar N_0^{C\beta}}.\tag 4.57$$
From Lemma 4.7 and the choice of $\sigma$, the expansion factor $C$ needs to 
satisfy 
$$\max(\frac{1}{\sigma}, 6(d+\nu))\ll C< \bar N_0^{\sigma/2}.\tag 4.58$$

From Theorem A and (4.19), $\bar N_0$ needs further to verify 
$$\bar N_0>\max(Q, M_0),\tag 4.59$$
where $M_0=M_0(d+\nu,\beta,\sigma)$ as in Lemma 4.2.
Fix $\bar N_0$ satisfying (4.59), $C$ satisfying (4.58). (4.57) determines $\epsilon_0$, $\delta_0$.

For the scales $\bar N_0\leq N\leq\bar N_0^C$,
$\bar N_0$ satisfying (4.59), the estimates 
$$\aligned 
&\Vert [T_N^\theta]^{-1} \Vert < e^{N^\beta}\qquad (0<\beta<1),\\
&|[T_N^\theta]^{-1}(\ell, \ell')|<e^{-\alpha|\ell-\ell'|}\qquad (\alpha \text{ as in }(H3)),\\
&|\ell-\ell'|>N/10, \endaligned$$ 
for $0<\epsilon\leq\epsilon_0$, $0<\delta\leq\delta_0$.
are obtained using Neumann series
by shifting in $\theta$ only.
So
$$
\text {mes\,} X_{N, i}=1,\quad \text {mes\,}\Omega_{N, i} =0, (\bar N_0\leq N\leq \bar N_0^C),
$$
where $X_{N,i}$, $\Omega_{N,i}$ are as defined in (4.51),
and (4.52) hold.

Let
$$\aligned X{\overset\text{\rm def }\to=}&
\bigcap_N \bigcap_{i\in[-3N^C, 3N^C]^{d+\nu}}\, X_{N, i}\\
=&\bigcap_N \bigcap_{j\in[-3N^C, 3N^C]^d}\, X_{N, j},\endaligned\tag 4.60$$
where the second equality follows from the remark after Lemma 4.5, $j$ being the $\Bbb Z^d$ coordinate of $i$; 
$$\Omega{\overset\text{\rm def }\to=}
\bigcup_N \bigcup_{i\in[-3N^C, 3N^C]^{d+\nu}}\, \Omega_{N, i},\tag 4.61$$
where $X_{N, j}$, $\Omega_{N, i}$ are as defined in (4.51).

On $X$, $DC_{A, c}\backslash\Omega$, Lemma 4.7 is available with the initial $\gamma=\alpha$ from
Lemma 4.2 for iteration to all scales. Estimating the measure of $X$ and $\Omega$ using (4.53, 4.54, 4.60, 4.61), using the measure estimates on the bad set in $\theta$ from Lemma 4.2 and 4.7, we obtain the Lemma by taking $p'\simeq\Cal O_d(1) C^3$.
$\hfill\square$
\bigskip

\noindent 
\head {\bf 5. Invertibility of $T(y_i)$, $Q$-equations and determination of $\omega$}
\endhead

Fix $x\in X\subset \Bbb R^{\Bbb Z^d}\backslash \Bbb R^\nu$, defined in (4.12), which generates a corresponding set
$\Omega$ as in Lemma 4.1. The main work of this section is to convert the measure estimates in $\theta$ for fixed
$\omega\in\Bbb R^{\nu}\backslash\Omega$, fixed $\Cal V\in \Bbb R^{\nu}$ in Lemma 4.1 into measure estimates in
$\omega=\omega(\Cal V)\in\Bbb R^\nu$ and extend to
$(\omega,\Cal V)\in \Bbb R^{2\nu}$ in the neighborhood of $\omega=\omega(\Cal V)$, while
keeping $\theta$ fixed:
$\theta=0$ and addressing the original family of linearized operators $T_N(y_i)$ defined in (4.2-4.5), where $y_i$
is the
$i^{\text{th}}$ approximant to the $P$-equations in (4.1).
This is possible because $T_N^\theta$ is only a function of $n\cdot\omega+\theta$, $T_N^\theta=
T'(n\cdot\omega+\theta)$, cf. (4.2-4.4, 4.7).

Since $\theta$ is now fixed at $0$, before doing the conversion, we need to make a further restriction in
$X$ in order that the spectrum of the various restricted random Schr\"odinger operators stay away from $0$.
This is needed when $n=0$ and we cannot vary $\omega$ to have invertibility of the linearized operators, cf. (4.3).

So we modify the definition of $X'_{N}$ in (4.23, 4.24) to also include the condition
$$\text{dist }(\sigma(\bar\Lambda(N')),0)>e^{-N_0^{\beta'}}\quad (0<\beta'<\beta)\tag 5.0$$
for all $\bar\Lambda(N')\cap[-N, N]^d\neq\emptyset$, ($N_0\leq N'\leq2N_0$), i.e., we require (4.30)
to hold also when $\lambda=\mu_j$, eigenvalues of restricted random Schr\"odinger operators and not
just differences of pairs of eigenvalues. In view of the Wegner estimate (A7), this leaves the measure
of the set $X'_N$ in (4.25) essentially unchanged.

This generates the restricted probability set $\tilde X\subset X\subset\Bbb R^{\Bbb Z^d}\backslash\Bbb R^\nu$,
on which Lemma 4.1 holds. Rename $\tilde X$ as $X$ and let 
$\Omega$ be its corresponding frequency set.
 
To do the conversion, we need to supplement the measure estimates in (4.14) by the fact that the bad sets $\Cal
B(N)$ defined from (4.9) has a semi-algebraic description in terms of $(\omega, \Cal V, \theta)$, enabling us to use
the decomposition Lemma 9.9 of [B5].
Once we have the necessary estimates on $[T_N (y_i)]^{-1}$ after removing a small set  in $\omega_i=\omega_i(\Cal V)$
(we put the suffix $i$ here to stress that it is the $i^{\text {th}}$ approximation), we
construct the next approximant $y_{i+1}$ according to (3.9), which in turn is used to construct
$\omega_{i+1}=\omega_{i+1}(\Cal V)$.

In this section, we primarily address the invertibility of $T_N(y_i), N=M^{i+1}$.
Since the estimate on $T_N(y_i)$ and the construction of $y_{i+1}$ are interconnected, it is good at
this point to lay down the complete induction hypothesis. (In section 4, we used first part of the induction 
hypothesis (H1, 3) and the first inequality of (H2)
to derive Lemma 4.1.) The first few approximations are obtained
by using direct perturbation series in $\epsilon$, $\delta$, (see (4.9), as well as the text and the remark
afterwards). So $0<\epsilon\leq\epsilon_0\ll 1$, $0<\delta\leq\delta_0\ll 1$, and we have $\alpha=\Cal O(1)|\log(\epsilon+\delta)|$.

On the {\it entire} ($\omega,\Cal V$) parameter space, we assume: 

\itemitem {(H1)} $\text{supp\,} y_i\subset [-M^i, M^i]^{d+\nu}$ ($i\geq 1$),

\itemitem {(H2)} $\Vert \Delta_i y\Vert <\delta_i$,\, $\Vert \partial \Delta_i y\Vert <\bar \delta_i$,

\noindent where $\Vert \ \Vert$ stands for $\sup_{\omega, \Cal V} \Vert \ \Vert_{\ell^2(\Bbb Z^{d+\nu})}$ 
(Recall that we identify $y$ with
$\hat y$; $\partial$ refers to derivation in $\omega$ or $\Cal V$.).
$\delta_i, \bar\delta_i$ will be shown to satisfy: 
$$\aligned \delta_i&<\sqrt{\epsilon+\delta}M^{-(\frac {4}{3})^{i}},\\
\bar\delta_i&<\sqrt{\epsilon+\delta}M^{-\frac{1}{2}(\frac {4}{3})^{i}}.\endaligned$$

\itemitem{(H3)} $|y_i(k)|< e^{-\alpha|k|}$ for some $\alpha>0$

\noindent (The constant in front of the exponential in (H3) is $1$, because $|a_\ell|\ll 1$, $\ell=1,...,\nu$, see
(1.9).)

From (H2), $y$ is a $C^1$ function on $(\omega,\Cal V)$.
Application of the implicit function theorem to the $Q$-equations in (2.12) with
$$
y=\pmatrix u\\ v\endpmatrix =\pmatrix \hat u\\ \hat v\endpmatrix =y_i\tag 5.1
$$
yields
$$
\omega_i=\Cal V+(\ve+\delta)\vp_i(\Cal V), \text { with } \Vert \partial\vp_i\Vert< C, \tag 5.2
$$
whose graph we denote by $\Gamma_i$. Recall that a priori, $y_i$ are only defined on certain intervals in 
$(\omega, \Cal V)$ space. It is in order to use the implicit function theorem that we extend $y_i$ to
the entire $(\omega, \Cal V)$ space, using the estimates on $\partial y_i$ in (H2).
(H2) and (2.12) imply, moreover, that
$$
\Vert\vp_i-\vp_{i-1} \Vert\leq \Cal O(1)\Vert y_i-y_{i-1}\Vert<\Cal O(1)
\delta_i, \tag 5.3
$$
which in turn implies that
$$
\Gamma_i \text { is an $(\epsilon+\delta)\delta_i$-approximation of $\Gamma_{i-1}$}.\tag 5.4
$$

At each stage $i$, define 
$$M_0=\Cal O(1)(i+1)^{C/2}(\log M)^{C/2}\tag 5.5$$
for some $C>0$ and $\Omega_{M_0, k}$ as defined in (4.51). Unlike  (H1-3), the following 
hypothesis is only assumed to hold on {\it certain} intervals in $\Bbb R^{2\nu}$, the
$(\omega, \Cal V)$-parameter space.

\itemitem{(H4)} There is a collection $\Sigma_i$ of intervals $I$ in $\Bbb R^{2\nu}$ of size $M^{-i^C}$,
same $C$ as in (5.5), such that

\itemitem{(H4, 0)} $(I\cap \Gamma_i)\subset(DC_{A, c}^{(i)}\backslash\Omega_i)$, where 
$$\aligned DC_{A, c}^{(i)}{\overset\text{\rm def }\to=}&DC_{A, c}(2N),\quad N=M^i,\\
\Omega_i{\overset\text{\rm def }\to=}&\cup_{k\in\Bbb Z^d\cap[-2M_0, 2M_0]^d}\Omega_{M_0, k},\quad i>i_0>0\endaligned
$$
(see the remark after (H4, iv) concerning $i_0$).

\itemitem
{(H4, i)} On $I\in \Sigma_i, y_i(\omega, \Cal V)$ is given by a rational function in ($\omega, \Cal V$) of degree at
most $M^{qi^3}(q\in \Bbb N^+)$.

\itemitem {(H4, ii)} For $(\omega, \Cal V)\in \bigcup_{I\in\Sigma_i} I$
$$
\align
&\Vert F(y_i)\Vert< \kappa_i\\
&\Vert\partial F(y_i)\Vert <\bar \kappa_i.\tag 5.6
\endalign
$$
In (6.20) $\kappa_i, \bar \kappa_i$ will be shown to satisfy:
$$\aligned \kappa_i&<\sqrt{\epsilon+\delta}M^{-(\frac {4}{3})^{i+2}},\\
\bar\kappa_i&<\sqrt{\epsilon+\delta}M^{-\frac{1}{2}(\frac {4}{3})^{i+2}}.\endaligned\tag 5.7$$

\itemitem {(H4, iii)} For $(\omega, \Cal V)\in \bigcup_{I\in\Sigma_i}I$,
$$
T_N=T_N(y_{i-1}), N=M^i
$$
satisfies
$$
\align
\Vert T_N^{-1}\Vert\leq M^{i^C}\\
|T_N^{-1}(k, k')|&\leq e^{-\alpha|k-k'|} \text { for } |k-k'|> i^C.
\endalign
$$

\itemitem
{(H4, iv)} Each $I\in \Sigma_i$ is contained in an interval $I'\in \Sigma_{i-1}$ and
$$
\text{mes}_\nu \big(\Gamma_i\bigcap (\bigcup_{I\in\Sigma_{i-1}} I\backslash \bigcup_{I\in\Sigma_i} I')\big)
< M^{-i/5}.
$$
\noindent {\it Remark.} (H4, 0) is only needed for $i>i_0$ to ensure the availability of Lemma 4.1. Up
to stage $i_0$, we use direct $\epsilon$, $\delta$ perturbation series, where the Diophantine property
of $\omega$ is not required, (cf. Lemma 4.2).

Contrary to related estimates on $ \Delta_i y$,\, $\partial\Delta_iy$ in (H2). (5.6, 5.7) {\it cannot} be extended
to the entire $(\omega, \Cal V)$ space. This is because, as mentioned  earlier, outside the intervals in
$\Sigma_i$,
$y_i$ are no longer close to solutions to the $P$-equations in (2.16).

\noindent
{\it Invertibility of} $T_N(y_i), N=M^{i+1}$

Assume (H1-4) hold at stage $i$. To construct $y_{i+1}$, we need to control
$$
[T_N (y_i)]^{-1}, N=M^{i+1},
$$
with a further restriction of the $(\omega, \Cal V)$-parameter set. This will give us (H4,iii) at stage $i+1$.

We accomplish this by covering $[-M^{i+1}, M^{i+1}]^{d+\nu}$ with $[-M^{i}, M^{i}]^{d+\nu}$ and intervals
$[-M_0, M_0]^{d+\nu}+k$, $M_0$ as in (5.5), $k\in\Bbb Z^{d+\nu}$, ${M^i}/2<|k|<M^{i+1}$ and using the resolvent
identity. We first estimate $[T_{M^i}(y_i)]^{-1}$. Fix $(\omega,\Cal V)\in\bigcup_{I\in\Sigma_i}I$. (H4,iii) at
stage $i$ gives
$$\aligned
\Vert [T_{M^i}(y_{i-1})]^{-1}\Vert&\leq M^{i^C},\\
|[T_{M^i}(y_{i-1})]^{-1}(k,k')|&\leq e^{-\alpha|k-k'|}\, (|k-k'|>i^C).\endaligned
\tag 5.8$$
We write
$$\aligned T_{M^i}(y_i)&=T_{M^i}(y_{i-1})+[T_{M^i}(y_i)-T_{M^i}(y_{i-1})]\\
&{\overset\text{def }\to=}A+B\endaligned\tag 5.9$$
From the first inequality in (5.8)
$$\Vert A^{-1}\Vert \leq M^{i^C}.\tag 5.10$$
The first inequality of (H2) at stage $i$ gives 
$$\Vert B\Vert <\Cal O(1) M^{-(\frac{4}{3})^i}.\tag 5.11$$
So $$\Vert [T_{M^i}(y_{i})]^{-1}\Vert\leq 2M^{i^C}, \text{ for } i>C^2.\tag 5.12$$

To obtain pointwise estimate on $T_{M^i}^{-1}(y_{i})$, we use (5.9) and resolvent series. $A^{-1}$
has off-diagonal decay from (5.8), $B$ has off-diagonal decay from (H3) at stage $i$. Iterating the
resolvent series and using (5.12), we obtain 
$$|[T_{M^i}(y_{i})]^{-1}(k,k')|\leq e^{-\alpha'|k-k'|}\, (|k-k'|>i^C),\tag 5.13$$
with $\alpha'=\alpha-M^{-i\delta'}>\alpha/2$ ($\delta'>0$),  uniformly in $i$.

We now study $|[T(y_{i})]^{-1}|$ on the $M_0$ intervals, $M_0$ as in (5.5).
We distinguish $2$ types of $M_0$ intervals $J$ in
$[-M^{i+1}, M^{i+1}]^{d+\nu}$:
\item{$\bullet$} $J\cap [-M_0, M_0]^d\times [-M^{i+1}, M^{i+1}]^\nu=\emptyset$ ($a$)
\item{$\bullet$} $J\cap [-M_0, M_0]^d\times [-M^{i+1}, M^{i+1}]^\nu\neq\emptyset$ ($b$)

\noindent For type $a$, we use direct perturbation in view of (H3). For type $b$, we use a more delicate
construction. We write $M_0=M^{\tilde i_0}$. ($M_0$ is chosen in order that the total degree of the semi-algebraic
set describing the bad set in $\omega$, $\Cal V$, $\theta$ is 
not too large.) 

The $M_0$
intervals are at
$[\tilde i_0]^{th}$ scale. On 
$M_0$ intervals of type $b$, we use Lemma 4.1, which is for $T_{M_0}(y_{i_0})$, $i_0\simeq\frac{\log M}{\log
b}\tilde i_0<i$ similar to (4.9). Using a decomposition lemma ([Lemma 9.9, B5] restated here as Lemma 5.3)
to make appropriate incisions in the $(\omega, \Cal V)$-parameter space,
applying (H2) between $y_{i_0}$ and $y_i$ and combining with estimates on type $a$ intervals, we obtain
$$\aligned
&\Vert (R_JT(y_i)R_J)^{-1}\Vert < e^{M_0^\beta}\qquad (0<\beta<1)\\
&|(R_JT(y_i)R_J)^{-1} (k, k')|<e^{-\alpha''|k-k'|}, k, k'\in J, |k-k'|>\frac{M_0}{10},\endaligned\tag 5.14$$
where $\alpha''=\alpha-M^{-i\delta''}$ ($\delta''>0$),
for all
$$
J=[-M_0, M_0]^{d+\nu}+k, \frac 12 M^i<|k|< M^{i+1}.\tag 5.15
$$
This is the content of Lemma 5.2.
We delay its precise statement and proof momentarily.
We first prove

\proclaim
{Lemma 5.1} Assume (5.12-5.15) hold. $M_0$ is as in (5.5). Then
$$
\align
&\Vert [T_{M^{i+1}} (y_i)]^{-1}\Vert < \Cal O(1) M^{(i+1)^C}\tag 5.16\\
&|[T_{M^{i+1}}(y_i)]^{-1}(k, k')|< e^{-\tilde\alpha|k-k'|} \text{ for } |k-k'|>(i+1)^C,
\tag 5.17
\endalign
$$
with $\tilde\alpha=\alpha-M^{-(i+1)\tilde\delta}$, $\tilde\delta>0$. 
\endproclaim

\noindent
{\bf Proof.} (5.16, 5.17) are exercises in resolvent identity or equivalently using Lemma B in the appendix.
We first prove (5.16).
(5.17) then follows by using (5.16) and another application of the resolvent identity.
Let
$$
B_i=[-M^i, M^i]^{d+\nu}, B_{i+1} =[-M^{i+1},M^{i+1}]^{d+\nu}.
$$

For any $k, \ell \in B_{i+1}$, we have (Assume $T^{-1}_{B_{i+1}}$ is defined.)
$$
\align
&T^{-1}_{B_{i+1}}(k, \ell)\\
&=T^{-1}_{W(k)} (k, \ell)+\sum_{\Sb k'\in \partial_*W(k)\\ |k''-k'|=1\\ k''\in B_{i+1}\backslash W(k)\endSb}
T^{-1}_{W(k)}(k, k')T^{-1}_{B_{i+1}}(k'',\ell)\tag 5.18
\endalign
$$
for $|k|\leq \frac 12 M^i, W(k)=B_i$,
for $|k|>\frac 12 M^i$, $W(k)$ is a size $M_0$ interval. 
It is easy to see that $\forall k'$, $\exists W(k)$, such that
$\text{dist}(k',\partial_* W(k))\geq M_0$, where $\partial_*W(k)$ is the interior boundary of $W(k)$, 
relative to $B_{i+1}$. Summing over $\ell \in B_{i+1}$ yields
$$
\align
&\sup_{k\in B_{i+1}}\sum_{\ell\in B_{i+1}}|T^{-1}_{B_{i+1}} (k, \ell)| \leq \sup_
{k\in B_{i+1}}\sum_{\ell\in W(k)}\Vert T^{-1}_{W(k)}\Vert\\
&+\sup_{k\in B_{i+1}} \sum_{k'\in \partial_*W(k)}|T^{-1}_{W(k)} (k, k')| \sup_{k''\in B_{i+1}} 
\sum_{\ell\in B_{i+1}}|T^{-1}_{B_{i+1}}(k'', \ell)|.
\tag 5.19
\endalign
$$

Using (5.12-5.15), we have
$$
\align
&\sup_{k\in B_{i+1}}\ \sum_{\ell\in B_{i+1}}|T^{-1}_{B_{i+1}}\, (k, \ell)|\leq 2M^{i^C}
\cdot (2M^i+1)^{d+\nu}\\
&+\Cal O(1) e^{-cM_0} M_0^{d+\nu-1}\sup_{k''\in B_{i+1}} \ \sum_{\ell\in B_{i+1}}|T^{-1}_{B_{i+1}}(k'',
\ell)|,
\tag 5.20
\endalign
$$
since $M^i\gg M_0$. So
$$
\Vert T^{-1}_{B_{i+1}}\Vert \leq M^{(i+1)^C},
$$
which is (5.16).

To obtain (5.17), we retrace our steps back to (5.18) and restrict to $k, \ell$, such that $|k-\ell|>(i+1)^C\gg
M_0$, in view of (5.5).
Iterating (5.18) along the path from $k$ to $\ell$ using $B_i$, $J$ and using (5.16) for the last factor, we
obtain (5.17). $\hfill\square$

(5.16, 5.17) will be the conclusion of (H4, iii) at stage $i+1$ once we specify the new set of intervals
$\Sigma_{i+1}$, on which they hold. As alluded to earlier, $\Sigma_{i+1}$ will be determined from
$\Sigma_{i}$ and the new restriction on $(\omega, \Cal V)$ in order that (5.14, 5.15) hold.

\noindent
{\it Determination of $\Sigma_{i+1}$}

\proclaim
{Lemma 5.2} Assume (H1-4) at stage $i$.
There exists $\tilde \Gamma_i\subset \Gamma_i$, $\text{\rm mes}_\nu \tilde \Gamma_i <M^{-i/4}$, such that
(5.14, 5.15) hold on
$$
\bigcup_{I\in\Lambda_i} (I\cap (\Gamma_i\backslash \tilde\Gamma_i).\tag 5.21
$$
\endproclaim

The proof of (5.21) relies on the measure estimates in Lemma 4.1 and semi-algebraic description of the bad set.
We need the following decomposition lemma, which is proven in [B5, Lemma 9.9].

\proclaim{Lemma 5.3}
Let $\Cal S\subset [0, 1]^{2n}$ be a semi-algebraic set of degree $B$ and $\text{\rm mes}_{2n} \Cal S <\eta, \log B\ll
\log 1/\eta$.
Denote by $(x, y)\in [0, 1]^n\times [0, 1]^n$ the product variable.
Fix $\ve>\eta^{1/2n}$.
Then there is a decomposition
$$
\Cal S =\Cal S_1 \bigcup \Cal S_2,
$$
with $\Cal S_1$ satisfying
$$
|\text{Proj}_x \Cal S_1|<B^K\epsilon \quad (K>0)\tag 5.22
$$
and $\Cal S_2$ satisfying the transversality property
$$
\text{\rm mes}_n(\Cal S_2\cap L)< B^K\epsilon^{-1} \eta^{1/2n}\quad (K>0),\tag 5.23
$$
for any $n$-dimensional hyperplane $L$ such that
$$
\max_{1\leq j\leq n}|\text{Proj}_L (e_j)|< \frac 1{100}\epsilon\tag 5.24
$$
where $e_j$ are the basis vectors for the $x$-coordinates.
\endproclaim

\noindent
{\bf Proof of Lemma 5.2.}
Assume (5.14, 5.15) hold with $T(y_{i_0})$ replacing $T(y_i)$, where as before
$$i_0\simeq\frac{\log M}{\log b}\tilde i_0\sim \frac{\log M_0}{\log b} <i,\tag 5.25$$
$$\aligned
\Vert (R_JT(y_{i_0})R_J)^{-1}\Vert < e^{M_0^\beta},\qquad (0<\beta<1)\\
|(R_JT(y_{i_0})R_J)^{-1} (k, k')|<e^{-\alpha|k-k'|}, k, k'\in J, |k-k'|>\frac{M_0}{10},\endaligned\tag 5.26
$$
for all
$$
J=[-M_0, M_0]^{d+\nu}+k, \frac 12 M^i<|k|< M^{i+1}.\tag 5.27
$$
Recall that we are at stage $i$, so (H2) is satisfied. Hence
$$\align \Vert T(y_i)-T(y_{i_0})\Vert &<\Cal O(1)\delta_{i_0}\\
&<e^{-\alpha M_0}\tag 5.28\endalign$$
using (5.25). This in turn implies that (5.14, 5.15) hold. So we only need to prove (5.26). ((5.28) is in
fact a reason for the choice of $M_0$ in (5.5).)

Fix $I\in\Sigma_{i_0}$. 
For type $a$ intervals, using (H3) for $y_{i_0}$, the off-diagonal elements in the $n$-direction of
$S$ has exponential decay and $\Vert S\Vert\leq\Cal O(\delta e^{-\alpha M_0})$. We use A. L. for the random
Schr\"odinger operator $\epsilon\Delta_j+V_j$ to obtain (5.26) as follows.

To obtain the first estimate in (5.26), we make direct incisions in the frequency space. Let $J$ be a type $a$
interval. We require
$$\aligned &|\pm n\cdot\omega+\epsilon\Delta_j+V_j|\\
=&|\pm n\cdot\omega+\mu_j|\\
\geq&2e^{-M_0^{\beta}},\endaligned$$
for all $(n,j)\in J$, where $\mu_j$ are the eigenvalues of $\epsilon\Delta_j+V_j$ restricted to the projection of
$J$ onto $\Bbb Z^d$.

When $n\neq 0$, this amounts to taking away a set in $\Omega$ of measure 
$\leq\Cal O(1) e^{-M_0^\beta}\cdot M_0^{d+\nu}$. When $n=0$, this is satisfied for $x\in X$ in view of (5.0).
So 
$$\Vert (R_JT(y_{i_0})R_j)^{-1}\Vert<e^{M_0^{\beta}},$$
by using the exponential estimates on $S$. This gives the first estimate in (5.26).

To obtain the second estimate in (5.26), we use Anderson localization, i.e., on the probability set $X$ defined in
(4.12), for all $E$ (here $E=n\cdot\omega$) there exists at most $1$ pairwise disjoint bad elementary regions
in $\Bbb Z^d$ of size $M_0^{1/C}$ in the projection of $J$ onto $\Bbb Z^d$. A resolvent series in the $j$ direction
coupled with a resolvent series in the $n$ direction using the above $2$ estimates and the decay property of $S$ 
gives the second estimate in (5.26), cf. proof of Lemma 4.3.

Hence there is a set $\bar\Gamma\subset \Gamma_{i_0}\cap I$,
$$\aligned \text {\rm mes\,}_\nu \bar\Gamma&<\Cal O(1)e^{-M_0^{\beta}}M^{(i+1)(d+\nu)}M_0^{d+\nu}\quad (0<\beta<1)\\
&<M^{-i}\endaligned\tag 5.29
$$
such that outside $\bar\Gamma$, (5.26) hold for all $J\subset[-M^{i+1}, M^{i+1}]^{d+\nu}$
satisfying $J\cap [-M_0, M_0]^d\times [-M^{i+1}, M^{i+1}]^\nu=\emptyset$. 

To prove (5.26) for type $b$ intervals $J$,
we use Lemma 4.1 at scale $\tilde i_0=[\frac{\log M_0}{\log M}]$ and the decomposition
Lemma 5.3. We illustrate this on the interval $J=[-M_0, M_0]^{d+\nu}$.
We consider the set
$$
\align
\Cal S= \{(\omega, \Cal V, \theta)\in I\times\Bbb R|
\, \text {(4.9) fail for $T^\theta_{[-M_0, M_0]^{d+\nu}}(y_{i_0})$, 
i.e., with $N$}&\text{ replaced by $M_0, y_i$},\\
& \text{ replaced by $y_{i_0}$}\},\tag 5.30
\endalign
$$
where $I\in\Sigma_{i_0}$ is the same fixed interval as earlier.
(Recall that $x\in X\subset\Bbb R^{\Bbb Z^d}\backslash\Bbb R^{\nu}$ is fixed.)
Let $T^\theta_{M_0}(y_{i_0})$ denote
$T^\theta_{[-M_0, M_0]^{d+\nu}}(y_{i_0})$.
Each matrix element of $T^\theta_{M_0}(y_{i_0})$ is a rational function in $\omega, \Cal V$ of degree 
at most $2(p+1)M^{q{i_0}^3}(q \in \Bbb N^+)$ and linear in $\theta$ (see (H4, i), (4.2-4.4, 4.7)). 
As before, the condition in (5.30) can be expressed in terms of determinants and hence polynomials in the matrix
elements of $T^\theta_{M_0}(y_{i_0})$.
This gives that $\Cal S$ is semi-algebraic of total degree at most $M_0^{C(d+\nu)} M^{q i_0^3}$ on $\Bbb
R^{2\nu+1}$.

We now localize to $(\omega, \Cal V)\in \Gamma_{i_0}$.
We consider the set
$$
\Cal S'{\overset\text{def }\to=} \Cal S \bigcap\{(I\cap \Gamma_{i_0}\cap\{DC_{A, c}(2M_0)\backslash\Omega_{M_0}\})\times\Bbb
R\}
\subset \Bbb R^{\nu+1},\tag 5.31
$$
where $\Omega_{M_0}$ is as in (4.54) of Lemma 4.7 with $N$ replaced by $M_0$,
$\Gamma_{i_0}$ is determined by the $Q$-equations in (2.12), which are polynomials in $(\omega, \Cal V)$ of degree at
most $2(p+1)M^{q i^3_0}(q \in \Bbb N^+)$. From Lemma 4.7, (4.54),
$DC_{A, c}(2M_0)\backslash\Omega_{M_0}$ is determined by $M_0^{C^2(d+\nu)}\sim M^{C^2(d+\nu)i_0}$ number of monomials of degree $1$.
So $\Cal S'$ is semi-algebraic of total degree at most $M^{(q+1)i_0^3}$ on $\Bbb R^{\nu+1}$.

By Lemma 4.1, each section $\Cal S'(\omega, \Cal V)=\Cal S'\big(\omega(\Cal V), \Cal V\big)$
is of measure at most $e^{-M^\sigma_0}$ in $\theta$ $(0<\sigma<1/2)$.
So
$$
\text{mes}_{\nu+1} \Cal S'\leq e^{-M_0^\sigma}\qquad (0<\sigma<1/2).\tag 5.32
$$
Our aim is to estimate for $k\in\Bbb Z^\nu$, $\frac 12 M^i<|k|\leq M^{i+1}$,
$$
\text{mes}_\nu\{\Cal V|\big(\Cal V,\Cal V+(\epsilon+\delta)\vp_{i_0}(\Cal V), k\cdot (\Cal
V+(\epsilon+\delta)\vp_{i_0}(\Cal V))\big)
\in \Cal S'\}.\tag 5.33
$$
Since $\Gamma_{i_0}$ is a $C^1$ function, $(I\cap \Gamma_{i_0})\times\Bbb R$ may be identified with an
interval in
$\Bbb R^{\nu+1}$, say $[0, 1]^\nu\times\Bbb R$, $\Cal S'$ defined in (5.31) is a subset of $(I\cap \Gamma_{i_0})\times
\Bbb R$, and therefore can be identified with a subset in $[0, 1]^\nu\times\Bbb R$.
For the purpose of the application of Lemma 5.3, we identify $[0, 1]^\nu\times \Bbb R$ with $[0, 1]^\nu\times [0,
1]^{\nu-1}\times\Bbb R$ and $\Cal S'$ with $\Cal S'\times [0, 1]^{\nu-1}$.
Since $T$ is restricted to the interval $[-M_0, M_0]^{d+\nu}$, we may further restrict the interval
in $[0, 1]^\nu\times\Bbb R$ to be 
$$
[0, 1]^\nu\times [-\Cal O(\sqrt {d+\nu})M_0, \Cal O(\sqrt {d+\nu})M_0]\simeq [0, 1]^{2\nu-1}
\times [-\Cal O(\sqrt {d+\nu})M_0, \Cal O(\sqrt {d+\nu})M_0].\tag 5.34
$$

We decompose $[-\Cal O(\sqrt {d+\nu})M_0, \Cal O(\sqrt {d+\nu})M_0]$ into intervals of length $1$ and
identify each of them with $[0, 1]$.
Applying the decomposition Lemma 5.3 to each of these intervals and taking the union, we obtain a subset 
$\Gamma'\subset\Gamma_{i_0}
\cap I$ of measure
$$
\align
\text{mes}_{\nu} \Gamma' &<(M^{(q+1)i_0^3})^C M^{-i}\cdot \Cal O(\sqrt {d+\nu})M_0\\
&< M^{-i/2},\tag 5.35
\endalign
$$
(Recall $\Cal S'$ is of degree $B\leq M^{q i_0^3}(q\in\Bbb N^+)$ and we take $\epsilon=200M^{-i}$ in (5.22).)
such that for all $k\in\Bbb Z^\nu, |k|>\frac 12 M^i$
$$
\align
&\text{mes}_\nu\{\Cal V|(\Cal V,\Cal V+(\epsilon+\delta)\vp_{i_0}(\Cal V), k\cdot \big(\Cal
V+(\ve+\delta)\vp_{i_0}(\Cal V)\big)\\
&\in\Cal S \cap \{\big((\Gamma_{i_0}\backslash \Gamma')\cap I\cap \{DC_{A,
c}(2M_0)\backslash\Omega_{M_0}\}\big)\times
\Bbb R\}\\ &<e^{-M_0^{\sigma/2}}\qquad (0<\sigma<1/2),\tag 5.36
\endalign
$$
where $M_0$ is as in (5.5). 

Same estimates as in (5.35, 5.36) hold when 
$T_{[-M_0, M_0]^{d+\nu}}$ is replaced by $T_{[-M_0, M_0]^{d+\nu}+\ell}$, $\ell\in \Bbb Z^d\cap[-2M_0, 2M_0]^d$.
Therefore, there is a set $\Gamma''\subset \Gamma_{i_0}\cap I$
$$
\align
\text{mes}_\nu\Gamma''&<\Cal O(1) M^{-i/2}+\Cal O(1) M^{(i+1)(d+\nu)} e^{-M_0^{\sigma/2}}\\
&<M^{-i/3}\tag 5.37
\endalign
$$
(in view of the choice of $M_0$ in (5.5) and $C\sigma\gg 1$ from (4.58)), such that outside $\Gamma''$, (4.9) hold
for all intervals
$J$ of the form $[-M_0, M_0]^{d+\nu}+\ell, \ell\in\Bbb Z^d\cap [-2M_0, 2M_0]^d$, $y_{i'}=y_{i_0}$
and
$\theta=k\cdot\omega, k\in\Bbb Z^\nu, \frac 12 M^i<|k|<M^{i+1}$. (The condition:
$\cap \{DC_{A, c}(2M_0)\backslash\Omega_{M_0}\}$ in (5.36) does not require additional incisions in 
the frequency space, as 
(H4, 0) holds starting at stage $i_0+1$.)

Combined with (5.29) and previous perturbation argument of replacing $y_i$ by $y_{i_0}$ in (5.28),
this implies that $\exists \Gamma'''$, $\text {\rm mes\,}_\nu\Gamma'''<M^{-i/3}$, such that outside $\Gamma'''$,
(5.14, 5.15) hold for all $k$, $\frac 12
M^i<|k|<M^{i+1}$, and a fixed $I\in\Sigma_{i_0}$.

Letting $I$ range over $\Sigma_{i_0}$, (There can be at most $\Cal O(1) M^{i_0^C}$ of these intervals, as the
$(\omega, \Cal V)$-parameter space can be restricted to, say $[0, 1]^{2\nu}$.) the total measure removed from
$\Gamma_{i_0}$ is at most $M^{i_0^C} \cdot M^{-i/3}<M^{-i/4}$.
Since $\Gamma_{i_0}$ and $\Gamma_i$ are at distance $<\delta_{i_0} <e^{-\alpha M_0}$ from
(5.3), we obtain a subset $\tilde\Gamma_i\subset\Gamma_i$, $\text{mes}_\nu \tilde\Gamma_i<M^{-i/4}$ such that
(5.14, 5.15) hold for all $k$, $M^i/2<|k|< M^{i+1}$ and on
$$
\bigcup_{I\in\Sigma_{i_0}}\big(I\cap(\Gamma_i\backslash\tilde\Gamma_i)\big)\tag 5.38
$$
and hence on
$$\bigcup_{I\in\Sigma_i} \big(I\cap (\Gamma_i\backslash \tilde\Gamma_i)\big)
\tag 5.39
$$
by (H4, iv).
This proves the Lemma.
$\hfill\square$

Lemma 5.1 then gives that on the  set in (5.39), (5.16, 5.17) hold.
Clearly by perturbation, (5.16, 5.17) remain valid on a $M^{-(i+1)^C}$ neighborhood of (5.39), (since 
$M^{-(i+1)^C}\ll e^{-\alpha M_0}\sim e^{-\alpha(i+1)^{C/2}}$ by the choice of $M_0$ in (5.5)),
which in turn
generates a collection $\Lambda_{i+1}$ of intervals in $\Bbb R^{2\nu}$ of size $M^{-(i+1)^C}$, such that for
$(\omega, \Cal V)\in I\in \Lambda_{i+1}$, (5.16, 5.17) hold.
So (H4, iii) hold at stage $i+1$ with $\tilde\alpha=\alpha-M^{-(i+1)\tilde\delta}$ ($\tilde\delta>0$), replacing $\alpha$.
Moreover we have
$$
\align
&\text{mes}_\nu\big(\bigcup_{I\in\Lambda_i}(I\cap\Gamma_i)\backslash \bigcup_{I'\in \Lambda_{i+1}}(I'\cap\Gamma_i)\big)
\\
&\leq \text{mes}_\nu \tilde\Gamma_i< M^{-i/4}\tag 5.40
\endalign
$$
which will imply (H4, iv) at stage $i+1$, once we construct $y_{i+1}$ and hence $\Gamma_{i+1}$ using (H4, iii) at
stage $i+1$.
\bigskip
\heading
6. Construction of $y_{i+1}$ and completion of the assemblage
\endheading

\noindent{\it Construction of $y_{i+1}$}

Let $N=M^{i+1}$, for $(\omega,\Cal V)\in\cup_{I\in\Lambda_{i+1}}I$, define
$$\Delta_{i+1}y=y_{i+1}-y_i{\overset\text{def }\to=}-(T_N(y_i))^{-1}F(y_i),\tag 6.1$$
(previously (3.9)). In view of (3.1, 2.16, H4,i), this implies that $\Delta_{i+1}y$ is a
rational function in $(\omega,\Cal V)$ of degree at most 
$$\aligned &\Cal O(1) N^{d+\nu}M^{q i^3}+(2p+1) M^{q i^3}\\
<&M^{q (i+1)^3}\qquad (q\in\Bbb N^+).\endaligned \tag 6.2$$
(Recall $p<M$.) So (H4, i) holds at stage $i+1$. (5.16), (H4,ii) give 
$$\Vert \Delta_{i+1}y\Vert <M^{(i+1)^C}\kappa_i=\delta_{i+1}\tag 6.3$$
and
$$\aligned \Vert \partial (\Delta_{i+1}y)\Vert & <\Vert \partial T_N^{-1}\Vert 
\Vert F(y_i)\Vert+\Vert T_N^{-1}\Vert \Vert \partial F(y_i)\Vert\\
&<\Vert T_N^{-1}\Vert ^2\Vert y_i\Vert_{C^1}\kappa_i+\Vert T_N^{-1}\Vert \bar\kappa_i\\
&<M^{2(i+1)^C}\bar \kappa_i\\
&=\bar\delta_{i+1},\endaligned \tag 6.4$$
where we also used (H2). 
Next we obtain point wise estimate on $\Delta_{i+1}y$. From (6.1)
$$|\Delta_{i+1}y(k)|\leq\sum_{|k'|\leq N}|T_N^{-1}(k,k')||F(y_i)(k')|.\tag 6.5$$
(2.16) gives 
$$\aligned |F(y_i)(k')|&\leq\Cal O(1)\sum_{k_1+\cdots+k_{2p+1}=k'}|y_i(k_1)|\cdots |y_i(k_{2p+1})|\\
&\leq (CM)^{CM}|k'|^{(d+\nu)M}e^{-\alpha|k'|},\endaligned\tag 6.6$$
(since $p<M$). Substituting (5.17, 6.6) into (6.5), we then obtain
$$\aligned |\Delta_{i+1}y(k)|&\leq (CM)^{CM}\{
\sum_{|k-k'|<i^C}M^{i^C}|k'|^{(d+\nu)M}e^{-\alpha|k'|}\\
&\qquad + \sum_{|k-k'|\geq i^C}|k'|^{(d+\nu)M}e^{-\tilde\alpha(|k'|+|k-k'|)}\}\\
&<C'M^{2i^C}|k|^{(d+\nu)M}e^{-\tilde\alpha|k|}.\endaligned\tag 6.7$$
Using (6.3) for $k$, such that $\log |k|\lesssim i$ and (6.7) otherwise, we obtain
$$|y_{i+1}(k)|<e^{(\log|k|)^{C'}-\tilde\alpha|k|}\leq e^{-\bar\alpha|k|},\tag 6.8$$
with $\bar\alpha=\alpha-M^{-(i+1)\bar\delta}$ for some $\bar\delta>0$ independent of $i$, where we used
the estimate on $\tilde\alpha$ just above (5.40).
This shows that (H3) is essetially preserved at stage $i+1$. Here we used the fact that $\alpha=\Cal O(1)|\log(\epsilon+\delta)|$ and $0<\epsilon\ll 1$, $0<\delta\ll 1$.

Since the intervals in $\Lambda_{i+1}$ are of size $M^{-(i+1)^C}$, we may extend $\Delta_{i+1}y$
to the entire $(\omega,\Cal V)$ parameter space as follows. For any $I\in\Lambda_{i+1}$, let
$\tilde I\subset I$ be such that $\text{dist}(\Lambda_{i+1}^c,\tilde I)\sim (1/3) M^{-(i+1)^C}$. Set
$\Delta y'_{i+1}=\Delta y_{i+1}$ on $I$, $\Delta y'_{i+1}=0$ on $\Lambda_{i+1}^c$. Define a $C^1$
function
$$\Xi_{i+1}=
\cases 1\qquad \text{on}\, \tilde I\\
0\qquad\text {on}\, \Lambda_{i+1}^c.\endcases
\tag 6.9
$$
Define $$\Delta\tilde y_{i+1}=\Xi_{i+1}\Delta y'_{i+1}.\tag 6.10$$
$\Delta\tilde y_{i+1}$ is defined on the whole $(\omega,\Cal V)$-parameter space and satisfies
$$\aligned \Vert \partial \Delta\tilde y_{i+1}\Vert
&<3M^{(i+1)^C}\delta_{i+1}+M^{2(i+1)^C}\bar\kappa_i\\
&=\bar\delta_{i+1},\endaligned\tag 6.11$$
where the second contribution comes from (6.4).
Renaming $\Delta \tilde y_i$, $\Delta y_i$ and letting $y_{i+1}=y_i+\Delta y_i$, we have thus
shown that (H1-3) remain valid at stage $i+1$ with $\bar\alpha$ replacing $\alpha$.

From $y_{i+1}$, the $Q$-equations in (2.12) define $\Gamma_{i+1}$ at most at a distance 
$\delta_{i+1}\simeq M^{-b^{i+1}}\ll M^{-i/4}$ from $\Gamma_i$. Clearly (5.40) implies 
$$\text {mes}_\nu(\Gamma_{i+1}\cap (\cup_{I\in\Lambda_i}I\backslash\cup_{I'\in\Lambda_{i+1}}I'))<M^{-i/4}
<M^{-(i+1)/5},\tag 6.12$$
which is (H4,iv) at stage $i+1$.

It remains to verify the properties of $F(y_{i+1})$ in (H4,ii), stage $i+1$. From the Taylor series in 
(3.10),
$$F(y_{i+1})=-[(T-T_N)[T_N(y_i)]^{-1}]F(y_i)+\Cal O(1)\Vert \Delta_{i+1} y\Vert^2,\qquad N=M^{i+1}.\tag 6.13$$
By construction, (H1), $\text{supp} y_i\subset [-M^i, M^i]^{d+\nu}$, (2.16) gives therefore, 
$$\aligned \text{supp} F(y_i)&\subset [-(2p+1)M^i, (2p+1)M^i]^{d+\nu}\\
&\subset [-M^{i+1}/10, M^{i+1}/10]^{d+\nu}\\
&=[-N/10, N/10]^{d+\nu}.\endaligned\tag 6.14$$
So 
$$\aligned F(y_{i+1})&=[R_{\Bbb Z^{d+\nu}\backslash B(0,N)}TT_N^{-1}R_{B(0,N/10)}]F(y_i)+\Cal O(1)\Vert
\Delta_{i+1}y\Vert^2,\\
\Vert F(y_{i+1})\Vert &\leq \Vert R_{\Bbb Z^{d+\nu}\backslash B(0,N)}TT_N^{-1}R_{B(0,N/10)}\Vert
\Vert F(y_i)\Vert +\Cal O(1)\Vert\Delta_{i+1}y\Vert^2,\endaligned \tag 6.15$$
where $B(0,N)=[-N,N]^{d+\nu}$, $R_{B(0,N)}$, $R_{B(0,N/10)}$ are characteristic functions. Thus 
$$\aligned \Vert F(y_{i+1})\Vert &\leq e^{-\alpha N/3}\kappa_i+\Cal O(1)\delta^2_{i+1}\\
&=\kappa_{i+1},\endaligned\tag 6.16$$
where we used (H2,3,4,ii). 

Similarly,
$$\aligned \Cal O(1)\Vert \partial F(y_{i+1})\Vert &\leq \Vert T_N^{-1}\Vert \Vert\partial T\Vert\Vert F(y_i)\Vert+
\Vert \partial T_N^{-1}\Vert \Vert F(y_i)\Vert\Vert T\Vert\\
&\qquad +\Vert R_{\Bbb Z^{d+\nu}\backslash B(0,N)}TT_N^{-1}R_{B(0,N/10)}\Vert \Vert \partial F(y_i)\Vert\\
&\qquad +\Vert \Delta_{i+1}y\Vert \Vert \partial \Delta_{i+1}y\Vert \\
&<M^{2(i+1)^C}\kappa_i+e^{-\alpha M^{i+1}/3}\bar\kappa_i+\delta_{i+1}\bar\delta_{i+1},\endaligned\tag 6.17$$
and we may take 
$$\bar\kappa_{i+1}=\Cal O(1)(M^{2(i+1)^C}\kappa_i+e^{-\alpha
M^{i+1}/3}\bar\kappa_i+\delta_{i+1}\bar\delta_{i+1}).\tag 6.18$$

Summarizing (6.3, 6.4, 6.16-6.18), we have 
$$\cases \delta_{i+1}=M^{(i+1)^C}\kappa_i\\
\bar\delta_{i+1}=M^{2(i+1)^C}\bar\kappa_i\\
\kappa_{i+1}=e^{-\alpha M^{i+1}/3}\kappa_i+\Cal O(1)\delta^2_{i+1}\\
\bar\kappa_{i+1}
=\Cal O(1)(M^{2(i+1)^C}\kappa_i+e^{-\alpha
M^{i+1}/3}\bar\kappa_i+\delta_{i+1}\bar\delta_{i+1}).\endcases
\tag 6.19
$$
We start from $\kappa_0$, $\bar\kappa_0=\Cal O(1)(\epsilon+\delta)$. For $\epsilon+\delta$ 
small enough, (6.19) is satisfied for $i\geq 1$, if 
$$\cases \delta_i<\sqrt {\epsilon+\delta} M^{-(\frac {4}{3})^i}, \qquad \kappa_i<\sqrt
{\epsilon+\delta}M^{-(\frac {4}{3})^{i+2}},\\
\bar\delta_i<\sqrt {\epsilon+\delta}M^{-\frac{1}{2}(\frac {4}{3})^i}, \qquad \bar\kappa_i<\sqrt
{\epsilon+\delta}M^{-\frac{1}{2}(\frac {4}{3})^{i+2}}.\endcases\tag 6.20$$

\noindent{\it (H4,0) and initial input for the induction}

To ensure (H4,0) at stage $i+1$, we make further incisions. (This is in order that Lemma 4.1 remains 
at our disposal at a later stage.) On $\Gamma_{i+1}$, we need to eliminate $\omega$ such that
$$\cases |n\cdot\omega_{i+1}+\lambda_{jj'}|\leq e^{-M_0^{\beta'}}\quad (n\sim\Cal O(1) M_0), \\
\Vert n\cdot\omega_{i+1}\Vert_{\Bbb T}\leq\frac {c}{|n|^A}\qquad (0<|n|<2M^{i+1}),\endcases
\tag 6.21$$
where $M_0$ is as in (5.5, H(4, 0)), $\beta'=\beta/{\Cal O(1)}$, $\Cal O(1)$ is the same expansion
factor as in Lemma 4.6 (denoted $C$ there), $\lambda_{jj'}=\mu_j-\mu_{j'}$, $\mu_j$, $\mu_{j'}$ 
are the eigenvalues of the random Schr\"odinger operator $\epsilon\Delta_j+V_j$ restricted to the 
myriad  elementary regions of size $M_0^{1/{\Cal O(1)}}$ (the same expansion factor $\Cal O(1)$) in $[-3M_0,
3M_0]^d$ (see Lemmas 4.1, 4.2, 4.6, proof of Lemma 4.2, (4.8, 4.9), the remark after (4.9, 4.26, 4.54) and the
definition of
$\Omega_i$ in (H4, 0)). There are at most
$\Cal O(1) M_0^{C'd}$ ($C'>0$) of such differences of eigenvalues.

In view of (5.2) at stage $i+1$, the first equation in (6.21) removes a set $\bar\Gamma_{i+1}\subset
\Gamma_{i+1}$, 
$$\aligned \text{mes}_\nu\bar\Gamma_{i+1}&\leq e^{-M_0^{\beta''}}\qquad (0<\beta''<\beta')\\
&\sim e^{-[(i+1)^{C/2}(\log M)^{C/2}]^{\beta''}}\\
&\ll M^{-\frac{i+1}{5}},\endaligned\tag 6.22$$
using (5.5) and choosing $$C>\frac{\Cal O(1)\log M}{\beta},$$
which is always possible.

Since 
$$\Vert \omega_{i+1}-\omega_i\Vert\leq\delta_i=M^{-(4/3)^i}\ll\frac{1}{M^{iA}}\tag 6.23$$
from (5.2, 5.3, 6.20), we only need to remove $\omega_{i+1}$ such that 
$$\Vert n\cdot\omega_{i+1}\Vert_{\Bbb T}\leq\frac {c}{|n|^A}$$
for $M^i\leq|n|\leq M^{i+1}$, which removes a set $\tilde\Gamma_{i+1}\subset\Gamma_{i+1}$,
$$\text{mes}_\nu\tilde\Gamma_{i+1}\leq \frac{\Cal O(1)}{M^{iA}}\ll M^{-\frac{i+1}{5}}\tag 6.24$$

Rename $\alpha$ as $\alpha_i$, $\bar\alpha$ as $\alpha_{i+1}$. From (6.8), 
$\alpha_{i+1}=\alpha_i-M^{-(i+1)\bar\delta}>\alpha/2$ uniformly in $i$. Combining (6.22, 6.24), we have (H4,0) at stage
 $i+1$ and that (6.12) is preserved. We have thus made
a complete induction from stage $i$ to $i+1$.
\bigskip

\heading 7. Proof of the Theorem
\endheading

The ``proof of the theorem" is now just a matter of juxtaposing sections 4, 5, 6 and recalling the sequence
of events. We recount the spine of the argument.

We use the modified Newton scheme in (3.9) to construct approximate solutions:
$$\Delta_{i+1}y=y_{i+1}-y_i=-(T_N(y_i))^{-1}F(y_i),\, N=M^{i+1},\tag 7.1$$ 
where $T_N$ is $T$ restricted to $[-N,N]^{d+\nu}$, $T$ and $F$ are as in (3.3-3.5, 3.1). Assume we have obtained
the first $i$ approximations $y_1,...,y_i$ on a set of intervals 
$\Bbb R^{2\nu}\supset\Lambda_1\supset\cdots\supset\Lambda_i$. To obtain $y_{i+1}$, we need to control $(T_N(y_i))^{-1}$
with a further restriction to the new set of intervals $\Lambda_{i+1}$ in $(\omega,\Cal V)$ space. This is accomplished as
follows.

To estimate $T_N(y_i)$, we cover $[-M^{i+1}, M^{i+1}]^{d+\nu}$ with the interval $[-M^{i}, M^{i}]^{d+\nu}=I$
and smaller intervals $J=[-M_0, M_0]^{d+\nu}+k,\, M^i/2<|k|<M^{i+1}$, $M_0\sim(\log N)^{C/2}$ as in (5.5). 
$T_I^{-1}$ is ``good" on $\Lambda_i$ by using perturbation theory. 
The $J$ intervals are divided into $2$ types as in section 5, according to their distances to the $\Bbb Z^d$
axis (see equations ($a$) and ($b$) between (5.13, 5.14)).
$T_J^{-1}$ of type
$a$ is easily obtained by using the manifest exponential decay properties of $y_i$ and a direct incision in 
the frequency space. The main task is to control $T_J^{-1}$ of type $b$,
which leads to further incisions in the frequency space, hence the new set of intervals $\Lambda_{i+1}$.

Since $|J|\ll|I|$, we may consider $T_J(y_{i_0})$ instead of $T_J(y_i)$ for some $i_0\ll i$ as in (5.25). 
We add a parameter $\theta$ to $T_J(y_{i_0})$ and estimate the measure of the set of $\theta$,
on the complement of which, $[T_J^{(\theta)}]^{-1}$ is ``good". This is Lemma 4.1. We then use the decomposition
Lemma 5.3 to transfer the estimate in $\theta$ into estimates in $\omega$, giving rise to the
new set of intervals $\Lambda_{i+1}$.

On $\Lambda_{i+1}$ we construct $y_{i+1}$ according to (7.1). Using the $Q$-equations in (2.12), we obtain
$\Gamma_{i+1}$. The first $i_0$ approximations are constructed by using direct $\epsilon$, $\delta$ series,
in order that Lemma 4.7 and hence Lemma 4.1 are available: $i_0\simeq \frac{1}{\beta}\log|\log(\epsilon+\delta)|$ from (4.57)
and the third expression in (4.9) after setting $N=\tilde N_0^C$ and determining $i$, hence $i_0$.
(6.20) gives the
rate of convergence of this Newton scheme and hence the Theorem.
$\hfill\square$
\bigskip
\heading
Appendix:
Localization results for random Schr\"odinger operators
\endheading

Random Schr\"odinger operator is the operator
$$
H=\epsilon \Delta +V \text { on } \ell^2 (\Bbb Z^d),
$$
where $\epsilon>0$ is a parameter, $\Delta (i, j) =1$ if $|i-j|=1$ and zero otherwise, $V=\{v_i\} _{i\in\Bbb Z^d}$ is a
family of independent identically distributed (iid) random variables with common probability distribution $g$.
The spectrum of $H$ is given by
$$
\align
\sigma (H) &=\sigma(\epsilon\Delta)+\sigma(V)\\
&=[-2\epsilon d, 2\epsilon d]+ \text {supp } g,\quad a.s.
\endalign
$$

We summarize below the known results on Anderson localization, which are relevant for the present 
construction (cf. \cite{DJLS, vDK, GB, GK, Mi, Si}). This is an expanded and more complete version of the
appendix in \cite{BW}.

For any $L\in\Bbb N$, let $\Lambda_L(i)$ denote any elementary region in $\Bbb Z^d$ with diameter
$2L$, center $i\in\Bbb Z^d$  as defined in (4.10, 4.11) with $\Bbb Z^d$ replacing $\Bbb Z^{d+\nu}$.
Let $H_{\Lambda_L(i)}$ be $H$ restricted to $\Lambda_L(i)$.
Let $m>0, E\in \Bbb R$.
$\Lambda_L(i)$ is $(m, E)$-regular (for a fixed $V$) if $E\not\in \sigma (H_{\Lambda_l(i)})$ and
$$
|G_{\Lambda_L(i)}(E; j, j')|\leq e^{-m|j-j'|}\tag A1
$$
for all $j, j'\in \Lambda_L(i), |j-j'|> L/4$. The following theorem is an immediate corollary
of the corresponding theorem in \cite{vDK} pertaining to cubes, by covering elementary regions with cubes and then an
application of the resolvent equation (cf. Lemma B). 
\proclaim
{Theorem A}
Let $I\subset\Bbb R$ be a bounded interval.
Suppose that for some $L_0>0$, we have
$$
 \text {{\rm Prob} \{for any $E\in I$ either $\Lambda_{L_0}(i)$ or $\Lambda_{L_0}(j)$ is
$(m_0, E)$-regular\} $\geq 1-\frac 1{L_0^{2p'}}$},\tag A2
$$ 
for some $p'>d, m_0>0$, and any $i, j\in\Bbb Z^d,  |i-j|> 2L_0$
$$
 \text{\rm Prob\,}\{\text{\rm dist\,} \big(E, \sigma(H_{\Lambda_L(0)})\big) < e^{-L^\beta}\}\leq 1/{L^{q'}}\tag A3
$$
for some $\beta$ and $q$, $0<\beta<1$ and $$q'> 4p'+6d\tag A4$$ all $E$ with
$$
\text{\rm dist\,} (E, I)\leq \frac 12 e^{-L^\beta},
$$
and all $L\geq L_0$.
Then there exists $\alpha, 1<\alpha <2$, such that if we set $L_{k+1}=[L^\alpha_k]+1$, $k= 0, 1,2 \ldots$ and pick $m,
0<m<m_0$, there is $Q<\infty$, such that if $L_0>Q$, we have that for any $k=0, 1,2 \ldots$
$$
\text{{\rm Prob} \{for any $E\in I$ either $\Lambda_{L_k(i)} $ or $\Lambda_{L_k(j)}$ is $(m, E)$ regular\} 
$\geq 1-\frac 1{L_k^{2p'}}$}\tag A5
$$
for any $i, j\in\Bbb Z^d$ with $|i-j|> 2L_k$.
\endproclaim

\noindent
{\it Remark.}
On the same probability subspace,
$$
\text {dist} \big(\sigma\big(H_{\Lambda_{L_k}(i)}\big), \sigma \big(H_{\Lambda_{L_k}(j)}\big)\big)>
e^{-L_k^\beta},\qquad
\beta>0\tag A6
$$
if $|i-j|> 2L_k$.
This is part of the ingredient of the proof of Theorem A.

Let $S\subset\Bbb Z^d$ be an (arbitrary) finite set. Let $H_S$ be $H$ restricted to $S$. If the probability 
distribution is absolutely continuous with a bounded density $\tilde g$, we have the following Wegner lemma:
$$\text{Prob }\{\text{dist } (E, \sigma(H_S))\leq\kappa\}\leq C\kappa|S|\Vert\tilde g\Vert_\infty,\quad
C>0,\,\kappa>0.\tag A7$$

(A2) is verified if $\epsilon$ is sufficiently small. (A3, 5) are provided by (A7), if $\Vert\tilde
g\Vert_\infty<\infty$. More precisely, fix $0<\beta<1$, choose $q'$ and hence $L_0$ sufficiently large, then there
exist
$\epsilon$ sufficiently small such that (A2, 3) are verified. We note from (A4) that the larger the $q'$, the larger
the
$p'$ could be. In view of (A5),
$q'$ can be chosen large if $L_0$ is large. So $p'$ can always be large enough by choosing $\epsilon$
small enough to suit the purpose of the construction in this paper (cf. Proof of Lemma 4.1).

Theorem A implies that for $0<\epsilon\ll 1$, $\Vert\tilde g\Vert_\infty<\infty$, $\sigma(H)$ has pure point 
spectrum almost surely. The pure point spectrum is dense. However it is simple \cite{Si}. Let $\psi_n$ ($n\in\Bbb
Z^d$) be the $n^{\text {th}}$ eigenfunction of $H$, then 
$$|\psi_n(j)|\leq C_{n,\omega}e^{-m'|j|}\quad (0<m'<m).$$
Further improvement of technology (see \cite{A, DJLS1,2, GB, GK}) give in fact that 
$$|\psi_n(j)|\leq C_{\omega}P_\omega (j_{n,\omega})e^{-m'|j-j_{n,\omega}|},\tag A8$$
where the centers $j_{n,\omega}$ satisfy $|j_{n,\omega}|\gtrsim n^{1/d}$, and $P_\omega$ is a polynomial,
which only depends on $\omega$.

\noindent{\it A resolvent estimate.}
\proclaim {Lemma B}
Suppose $\Lambda\subset \Bbb Z^{d+\nu}$ is an arbitrary set with the following property: for every $x\in\Lambda$, there is a subset $W(x)\subset\Lambda$ with $x\in W(x)$, $\text {diam }(W(x))\leq N$ and 
such that Green's function $G_{W(x)}(E)$ satisfies for certain $t$, $N$, $A>0$
$$\align \Vert G_{W(x)}(E)\Vert &<A\tag B1\\
|G_{W(x)}(E;x,y)|&<e^{-tN}\quad\text {for all } y\in\partial_*W(x).\tag B2\endalign$$
Here $\partial_*W(x)$ is the interior boundary of $W(x)$ relative to $\Lambda$ given by 
$$\partial_*W(x)=\{y'\in W(x)|\exists z\in\Lambda\backslash W(x),\quad |z-y'|=1\}.\tag B3$$
Then 
$$\Vert G_{\Lambda}(E)\Vert<2N^2A$$
provided $4N^2e^{-tN}\leq\frac{1}{2}$.
\endproclaim
See \cite{BGS}, where it is stated as Lemma 2.2, for a proof using the resolvent equation. See also the
proof of Lemma 5.1 in section 5 of the present paper for an essentially identical exercise in resolvent equation.
\hfill $\square$

\enddocument
\end